\pdfoutput=1
\RequirePackage{ifpdf}
\ifpdf 
\documentclass[pdftex]{sigma}
\else
\documentclass{sigma}
\fi

\usepackage{xcolor}

\usepackage{nicefrac, faktor}
\usepackage[shortlabels]{enumitem}
\usepackage{mathrsfs}

\numberwithin{equation}{section}

\newtheorem{Theorem}{Theorem}[section]
\newtheorem*{Theorem*}{Theorem}
\newtheorem{Corollary}[Theorem]{Corollary}
\newtheorem{Lemma}[Theorem]{Lemma}
\newtheorem{Proposition}[Theorem]{Proposition}
\newtheorem{Assertion}[Theorem]{Assertion}
\newtheorem{Conjecture}[Theorem]{Conjecture}
 { \theoremstyle{definition}
\newtheorem{Definition}[Theorem]{Definition}

\newtheorem{Example}[Theorem]{Example}
\newtheorem{Remark}[Theorem]{Remark}
}

\newcommand{\F}{\mathbb{F}}
\newcommand{\Z}{\mathbb{Z}}
\newcommand{\R}{\mathbb{R}}
\newcommand{\C}{\mathbb{C}}
\newcommand{\N}{\mathbb{N}}
\newcommand{\SSS}{\mathbb{S}}
\newcommand{\TTT}{\mathbb{T}}

\newcommand{\SO}{\mathrm{SO}}
\newcommand{\charac}{\mathrm{char}}
\newcommand{\eps}{\varepsilon}
\newcommand{\Sc}{{\rm Sc}}

\newcommand{\SU}{\rm SU}
\newcommand{\id}{{\rm id}}
\newcommand{\MSO}{{\rm MSO}}
\newcommand{\sB}{{\sf B}}
\newcommand{\Tor}{{\rm Tor}}
\newcommand{\im}{{\rm im}}

\newcommand{\myicon}{$\,\,\,\triangleright$}

\newcommand{\no}{{\color{brown!70!black}$\nexists$}{\sf PSC}}
\newcommand{\yes}{{\color {blue!70!black}$\exists$}{\sf PSC}}

\newcommand{\Ho}{\mathrm{H}}

\begin{document}


\renewcommand{\thefootnote}{}

\newcommand{\arXivNumber}{2112.04825}

\renewcommand{\PaperNumber}{069}

\FirstPageHeading

\ShortArticleName{Torsion Obstructions to Positive Scalar Curvature}

\ArticleName{Torsion Obstructions to Positive Scalar Curvature\footnote{This paper is a~contribution to the Special Issue on Global Analysis on Manifolds in honor of Christian B\"ar for his 60th birthday. The~full collection is available at \href{https://www.emis.de/journals/SIGMA/Baer.html}{https://www.emis.de/journals/SIGMA/Baer.html}}}

\Author{Misha GROMOV~$^{\rm ab}$ and Bernhard HANKE~$^{\rm c}$}

\AuthorNameForHeading{M.~Gromov and B.~Hanke}

\Address{$^{\rm a)}$~Institut des Hautes \'Etudes Scientifiques, 91893 Bures-sur-Yvette, France}
\EmailD{\href{mailto:gromov@ihes.fr}{gromov@ihes.fr}}
\URLaddressD{\url{https://www.ihes.fr/~gromov/}}

\Address{$^{\rm b)}$~Courant Institute of Mathematical Sciences, New York University,\\
\hphantom{$^{\rm b)}$}~New York, NY~10012-1185, USA}

\Address{$^{\rm c)}$~Institut f\"ur Mathematik, University of Augsburg, 86135 Augsburg, Germany}
\EmailD{\href{hanke@math.uni-augsburg.de}{hanke@math.uni-augsburg.de}}
\URLaddressD{\url{https://www.math.uni-augsburg.de/diff/hanke/}}

\ArticleDates{Received November 07, 2023, in final form July 17, 2024; Published online July 30, 2024}

\Abstract{We study obstructions to the existence of Riemannian metrics of positive scalar curvature on closed smooth manifolds arising from torsion classes in the integral homology of their fundamental groups. As an application, we construct new examples of manifolds which do not admit positive scalar curvature metrics, but whose Cartesian products admit such metrics.}


\Classification{53C21; 53C23; 20J06; 53C12; 55N10}

\begin{flushright}
\begin{minipage}{60mm}
\it Dedicated to Christian B\"ar\\
 on the occasion of his 60th birthday
 \end{minipage}
\end{flushright}

\renewcommand{\thefootnote}{\arabic{footnote}}
\setcounter{footnote}{0}


\section{Introduction}

In his construction of counterexamples to the Gromov--Lawson--Rosenberg conjecture, Thomas Schick in \cite{Sch98} discovered a purely homological obstruction to the existence of metrics with positive scalar curvature on closed oriented smooth manifolds.
We say that a smooth manifold is~\yes, if it admits a Riemannian metric of positive scalar curvature, and \no, if it does not admit such a metric.
Schick's construction is based on the observation that by the descent method of Schoen--Yau \cite{S-Y1979,SY}, the \yes-property of a closed oriented manifold $V$ of dimension~${n \leq 7}$ implies that the Poincar\'e dual of a class $c \in H^{n-\ell}(V; \Z)$ which is the product of $1$-dimensional integral cohomology classes can be represented by a~closed oriented smoothly embedded \yes-submanifold $\Sigma \subset V$ of dimension $\ell$.
Now, if $\ell = 2$ and the fundamental class of $\Sigma$ can be detected by a product of two $1$-dimensional cohomology classes in $V$ with coefficients in some field $\F$, then $\Sigma$ is \no~by the Gauss--Bonnet theorem.
This cohomological condition on the fundamental class of $\Sigma$ is satisfied if there exist classes $c_1, c_2 \in H^1(V; \F)$ such that $c \cup c_1 \cup c_2 \in H^n(V; \F)$ evaluates to nonzero on the fundamental class of $V$.
It follows by contraposition that in this case $V$ is \no.

In the present paper, we generalize these ``toral'' \yes-obstructions in two directions.
First, we combine the torality condition on the Poincar\'e dual of $c$ with an enlargeability condition on $V$, see Definition \ref{def:enlargeable}.
The resulting largeness property depends only on the image of the fundamental class of $V$ in $H_n( \pi_1(V) ; \F)$ under the classifying map, see Proposition \ref{homol_inv}.
Second, we allow the case $\dim \Sigma = 3$ by exploiting the restricted nature of \yes-$3$-manifolds, see Proposition \ref{prop:toralthree}.
In Theorem \ref{thm:tornexist}, we obtain a general \yes-obstruction in this setting.

As an application, in Theorem \ref{thm:plenty}, we construct new examples of \no-manifolds whose Cartesian products are \yes, thus challenging our geometric intuition that \no-manifolds are ``large'' while \yes-manifolds are ``small''.

Some fundamental questions remain unanswered, most notably whether there exist toral \yes-manifolds of dimension $n \geq 4$ and with fundamental group $(\Z/p)^n$ for some odd prime $p$, see Conjecture \ref{main_conj} and \cite[Problem 9]{MP2021}.
We hope that this note will stimulate some interest in this and related problems, some of which we address in Section \ref{fourconjectures}.

In an appendix, we prove homological invariance of the \yes-property of closed oriented manifolds of dimension at least $5$ with non-spin universal covers, filling a gap in an earlier treatment in the literature \cite{RS}.
This result is of independent interest.

\section[Torality and tor-enlargeability]{$\boldsymbol{\F}$-torality and tor-$\boldsymbol{\ell/\F}$-enlargeability}

\begin{Definition} \label{deftoral}
Let $X$ be a topological space, let $\ell$ be a non-negative integer and let $\mathbb{F}$ be a~field, e.g., $\mathbb{Q}$ or the finite field $\mathbb F_{p}$ for some prime $p$.

In the line of \cite{Han20,Jo2004}, a homology class $h \in H_\ell(X; \F)$ is called {\it $\mathbb F$-toral} if it is detected by a~product of $1$-dimensional cohomology classes, i.e., there exist classes $c_1, \ldots, c_\ell \in H^1(X ; \mathbb F)$ such that
\[
 \langle c_1\cup \cdots \cup c_\ell, h \rangle \neq 0 \in \F.
\]
Here we use the Kronecker pairing (evaluation pairing) of cohomology and homology.
We agree that all $h \neq 0 \in H_0(X; \F)$ are $\F$-toral.
Otherwise, the class $h$ is called {\it $\F$-atoral}.
Obviously, the subset of $\F$-atoral classes forms a linear subspace of $H_{\ell}(X; \F)$.

A closed, i.e., compact-without-boundary, $\F$-orientable\footnote{This condition is empty for $\charac(\F) = 2$ and equivalent to orientability for $\charac(\F) \neq 2$.} $n$-dimensional manifold $V$ is called {\it $\F$-toral}, if, after choosing some $\F$-orientation, its fundamental class $[V]_{\F} \in H_{n}(V; \F)$ is $\F$-toral.
Otherwise, $V$ is called {\it $\F$-atoral.}
\end{Definition}

In the following, let $\SSS^n$ denote the unit $n$-sphere equipped with the standard round metric, and let $\infty \in \SSS^n$ be some base point.
We call a continuous map $f \colon X \to \SSS^n$ {\em compactly supported}, if there exists a compact subset $K \subset X$ with $f(X \setminus K) = \{ \infty\}$.

\begin{Remark} \label{rem:bemerkung}\quad
\begin{enumerate}[(a)]\itemsep=0pt
\item If $\charac(\F) \neq 2$, then the classes $c_1, \ldots, c_\ell$ appearing in Definition \ref{deftoral} span an exterior algebra on $\ell$ generators of dimension $1$ in the $\F$-algebra $H^*(X ; \F)$.
\item \label{overtorical} A closed oriented $n$-manifold $V$ is $\mathbb{Q}$-toral, if and only if it admits a map of non-zero degree to the $n$-torus $\TTT^n$.
In other words, $V$ is {\em overtorical}.
This follows from the universal coefficient formula $H^*(V; \mathbb{Q}) = H^*(V; \mathbb{Z}) \otimes \mathbb{Q}$ and the one-to-one correspondence of classes in $H^1(V; \mathbb{Z})$ and homotopy classes of maps $V \to S^1$.
\end{enumerate}
\end{Remark}

\begin{Definition} \label{def:enlargeable} Let $n \geq 1$ and $0 \leq \ell \leq n$.
An $\F$-orientable Riemannian $n$-manifold $(V,g)$ is called {\it tor-$\ell/\mathbb F$-enlargeable}, if for all $\varepsilon > 0$, there exists a Riemannian covering
$
 \bar V = \bar V_\varepsilon \to V
$
and a compactly supported smooth $\varepsilon$-Lipschitz map
$
 f = f_\varepsilon \colon \bar V \to \SSS^{n-\ell}
$,
which is of {\em non-zero $\F$-toral degree}.
By definition, this means that for some, hence every, regular value ${s \in \SSS^{n - \ell} \setminus \{\infty\}}$ of~$f$, the homology class $\big[f^{-1}(s)\big] \in H_{\ell}\big(\bar V; \F\big)$ represented by the closed $\F$-orientable submanifold~${f^{-1}(s) \subset \bar V}$ is $\F$-toral.

For $\ell = n$, this means that a finite cover of $V$ is $\F$-toral.
\end{Definition}

\begin{Remark}\quad
\begin{enumerate}[(a)]\itemsep=0pt
\item In the case of closed $V$, this property is independent of the choice of $g$.
Thus, it is an invariant of the diffeomorphism type of $V$.
\item In fact, if $V$ is closed, connected and $\F$-oriented, this property depends only on the image of the fundamental class $[V]_\F \in H_n(V; \F)$ under the classifying map $V \to {\sf B}\pi_1(V)$, see Proposition \ref{homol_inv}.
In particular, for closed $V$, tor-$\ell/\F$-enlargeability is an invariant of the homotopy type of $V$.
\item Tor-$0/ \mathbb{Q}$-enlargeable $n$-manifolds $V$ are {\em enlargeable} in the classical sense of \cite[Section 5]{G-L1983}, i.e., for all $\eps > 0$, there exists a Riemannian covering $\bar V \to V$ together with a compactly supported $\eps$-Lipschitz map $\bar V \to \SSS^n$ of non-zero degree.
\item ``Tor-$\ell/\F$-enlargeable'' generalizes ``$\mathcal{SYS}$-enlargeable'' from \cite{Gr2018}, where $\mathcal{SYS}$ stands for \textbf Schoen--\textbf Yau--\textbf Schick.
\end{enumerate}
\end{Remark}

\begin{Example} \label{ex:toral}\quad
\begin{enumerate}[(a)]\itemsep=0pt
\item \label{eins} For all $n \geq 1$ and all $\F$, the $n$-torus $V = \TTT^n = S^1 \times \cdots \times S^1$ is $\F$-toral.
\item \label{projspace} For all $n \geq 1$, the real projective space $\mathbb{R} P^n$ is $\F_2$-toral.
\item \label{ex:overtorical} Let $0 \leq \ell \leq n$ and let $V$ be a closed $\F$-orientable $n$-manifold together with a map $f \colon V \to {\TTT}^{n-\ell}$ of non-zero $\F$-toral degree.
Then $V$ is $\F$-toral and tor-$\ell/\F$-enlargeable.
Such $V$ are called {\em tor-$\ell/\F$-overtorical}.
\item Let $V$ be tor-$\ell / \F$-enlargeable and let $W$ be a an $\F$-orientable Riemannian manifold admitting a proper $1$-Lipschitz map $W \to V$ of non-zero $\F$-toral degree.
Then $W$ is tor-$\ell/\F$-enlargeable.
\end{enumerate}
\end{Example}

For a group $G$, let ${\sf B}G$ be the classifying space of $G$ and let $H_*(G) = H_*({\sf B}G)$ be the group homology of $G$.
For a connected topological space $X$, let
$\varphi \colon  X \to {\sf B} \pi_1(X)$
be the classifying map of the universal cover of~$X$.
It is unique up to homotopy.

\begin{Definition} \label{def:essential} Let $V$ be a closed connected manifold of dimension $n$.
\begin{enumerate}[(a)]\itemsep=0pt
\item Let $\F$ be a field.
$V$ is called {\em $\F$-essential} if it is $\F$-orientable and the classifying map $\varphi \colon V \to {\sf B}\pi_1(V)$ sends some fundamental class $[V] _\F \in H_{n}(V; \F)$ to a non-zero homology class in $H_n(\pi_1(V); \F)$.
\item $V$ is called {\em $($integrally$)$ essential} if it is orientable and $\varphi$ sends some integral fundamental class $[V] \in H_n(V;\Z)$ to a non-zero homology class in $H_n(\pi_1(V); \Z)$.
\end{enumerate}
\end{Definition}

Any oriented $\F$-essential manifold $V^n$ is essential, since $\varphi_*([V]_\F) \in H_n( \pi_1(V) ; \F)$ is just the $\F$-reduction of $\varphi_*([V]) \in H_n(\pi_1(V) ; \Z)$.
However, the last class may be a {\em torsion class}, see Example \ref{ex:torsion}.
Hence, the title of our paper.

In the following, for $n \geq 0$, let $\sigma^n \in H^{n}_{\rm c}(\SSS^{n} \setminus \{ \infty\}; \F)$ be the non-zero compactly supported cohomology class Poincar\'e dual to the class in $H_0( \SSS^{n} \setminus \{ \infty\}; \F) = \F$ represented by a point.
Note that for each compactly supported map $f \colon X \to \SSS^{n}$, the cohomology class $f^*(\sigma^n) \in H^{n}_{\rm c}(X; \F)$ is compactly supported.

\begin{Proposition} \label{special} Each connected tor-$\ell / \F$-enlargeable manifold $V$ is $\F$-essential.
\end{Proposition}

\begin{proof}
We construct a CW model of ${\sf B}\pi_1(V)$ by attaching cells of dimension $\geq 3$ to $V$, killing all homotopy groups in dimension $\geq 2$.
In this description, the classifying map $\varphi \colon V \to {\sf B}\pi_1(V)$ becomes an inclusion.

Let $n = \dim V$ and assume that $\varphi_*([V]_\F) = 0 \in H_n ( \pi_1(V) ; \F)$.
Then there exists a connected finite subcomplex $V \subset S \subset {\sf B}\pi_1(V)$ such that the inclusion map $\iota \colon V \hookrightarrow S$ satisfies $\iota_*([V]_{\F}) =0 \in H_n(S; \F)$.
Since $V$ and $S$ have the same $2$-skeleton, the map $\iota$ induces an isomorphism of fundamental groups.

By assumption and Poincar\'e duality, for all $\varepsilon > 0$, there exists a connected covering $\pi_V \colon \bar V = \bar V_\eps \to V$, a compactly supported smooth $\eps$-Lipschitz map $f_\eps \colon \bar V \to \SSS^{n-\ell}$, and classes $c_1, \ldots , c_{\ell} \in H^1\big(\bar V; \F\big)$ so that
\begin{equation} \label{thetadeg}
 \big( f_\eps^*\big(\sigma^{n-\ell}\big) \cup (c_1 \cup \cdots \cup c_\ell ) \big) \cap \pi_V^!( [V]_{\F}) \neq 0 \in H_0 \big( \bar V; \F\big) .
\end{equation}
 Here $H_*^{\rm lf}$ is locally finite homology and $\pi_V^{!} \colon H_n(V ; \F) \to H^{\rm lf}_n\big( \bar V; \F\big)$ is the homological transfer.
Note that $f_\eps^*\big(\sigma^{n-\ell}\big) \cup (c_1 \cup \cdots \cup c_\ell ) \in H_c^n\big( \bar V; \F\big)$, so that the class considered in~\eqref{thetadeg} indeed lies in $H_0\big(\bar V; \F\big)$.

 Let $G := \pi_1\big(\bar V\big) < \pi_1(V)$ and let $\pi \colon {\sf B}G \to {\sf B}\pi_1(V)$ be the corresponding covering map.
Let~${\bar S := \pi^{-1}(S) \subset {\sf B}G}$.
Since $S \hookrightarrow {\sf B}\pi_1(V)$ induces an isomorphism of fundamental groups, $\bar S$ is connected and the inclusion $\bar V \hookrightarrow \bar S$ induces an isomorphism of fundamental groups.
Hence, we obtain $H^1\big(\bar S; \F\big) \cong H^1\big(\bar V; \F\big)$, and the classes $c_1, \ldots, c_\ell \in H^1\big(\bar V; \F\big)$ extend to classes $c_1', \ldots, c_\ell' \in H^1\big(\bar S; \F\big)$.

By \cite[Lemma 3.2]{BH10} and since $S$ is a finite complex, one can choose $\eps > 0$ at the beginning of this argument so small that the map $f_\eps \colon \bar V \to \SSS^{n-\ell}$ extends to a compactly supported map $F \colon \bar S \to \SSS^{n-\ell}$.

Let $\pi_S^! \colon H_n(S; \F) \to H^{\rm lf}_n\big(\bar S; \F\big)$ be the homological transfer map for the covering $\pi_S \colon \bar S \to S$.
Then \eqref{thetadeg} together with the naturality of the transfer map implies
\[
 \big( F^*\big(\sigma^{n-\ell}\big) \cup ( c_1' \cup \cdots \cup c_\ell') \big) \cap \pi_S^!(\iota_*( [V]_{\F})) \neq 0 \in H_0\big(\bar S; \F\big) .
\]
This contradicts our assumption $\iota_*([V]_\F) = 0 \in H_n(S; \F)$.
\end{proof}

Extending this line of thought, in Proposition \ref{homol_inv}, we will show a homological invariance property of tor-$\ell / \F$-enlargability.
This is similar in spirit to \cite[Section 3]{BH10}.
In the following, we equip connected simplicial complexes with their canonical path metrics whose restriction to each $k$-simplex is the standard metric on $\Delta^k \subset \R^{k+1}$.

\begin{Definition} \label{largehom}
Let $0 \leq \ell \leq n$.
Let $X$ be a connected simplicial complex with finitely presented fundamental group.
A homology class $h \in H_n(X; \F)$ is called {\em tor-$\ell / \F$-enleargeable}, if there exists a finite connected subcomplex $S \subset X$ and a class $h_S \in H_n(S; \F)$ such that
\begin{enumerate}[(a)]\itemsep=0pt
\item \label{isofund} The inclusion $S \hookrightarrow X$ induces an isomorphism of fundamental groups\footnote{We implicitly assume that $S$ contains the base point of $X$.} and the induced map $H_n(S ; \F) \to H_n(X; \F)$ sends $h_S$ to $h$;
\item \label{shrinkingmap} for all $\eps > 0$, there exists a connected cover $\pi_S \colon \bar S = \bar S_\eps \to S$, a compactly supported $\eps$-Lipschitz map $f_\eps \colon \bar S \to \SSS^{n-\ell}$ and classes $c_1 , \ldots, c_\ell \in H^1\big(\bar S; \F\big)$ that satisfy
 \begin{equation*}
 \big( f_\eps^*\big(\sigma^{n-\ell}\big) \cup (c_1 \cup \cdots \cup c_\ell) \big) \cap \pi_S^!(h_S) \neq 0 \in H_{0}\big(\bar S; \F\big) .
\end{equation*}
\end{enumerate}
Let
\smash{$
 H_*^{\text{\rm tor-}\ell/\F\text{\rm-enl}}(X; \F) \subset H_*(X ; \F)
$}
be the subset of tor-$\ell/\F$-enlargeable homology classes.
\end{Definition}

The complex $S$ in Definition \ref{largehom} can be made arbitrarily large:

\begin{Lemma} \label{subcomplexinv} Let $h\in H_n(X; \F)$ be tor-$\ell / \F$-enlargeable.
Let $S \subset X$ and $h_S \in H_n(S;\F)$ as in Definition {\rm\ref{largehom}} and let $K \subset X$ be a finite subcomplex.

Then there exists a finite connected subcomplex $S \cup K \subset T \subset X$ such that the inclusion~${T \hookrightarrow X}$ induces an isomorphism of fundamental groups.

Given such a subcomplex $T \subset X$, let $h_T \in H_n(T; \F)$ be the image of $h_S$ under the inclusion~${S \hookrightarrow T}$.
Then $T$ and $h_T$ satisfy properties {\rm\ref{isofund}} and {\rm\ref{shrinkingmap}} in Definition {\rm\ref{largehom}}.
\end{Lemma}

\begin{proof} The existence of $T$ is clear.
Obviously, the induced map $H_n(T; \F) \to H_n(X; \F)$ sends~$h_T$ to $h$.
Let $\eps > 0$, let $f_{\eps} \colon \bar S \to \SSS^{n-\ell}$ be an $\eps$-contracting compactly supported map where $\bar S \to S$ is a connected cover and let $c_1, \ldots, c_\ell \in H^1\big(\bar S; \F\big)$ such that
\begin{equation} \label{nonzero}
 \big( f_\eps^*\big(\sigma^{n-\ell}\big) \cup (c_1 \cup \cdots \cup c_\ell) \big) \cap \pi_S^!(h_S) \neq 0 \in H_{0}\big(\bar S; \F\big) .
\end{equation}
Let $\bar T \to T$ be the uniquely determined connected cover which restricts to $\bar S \to S$.
The inclusion~${\iota\colon \bar S \hookrightarrow \bar T}$ induces an isomorphism of fundamental groups, hence an isomorphism $H^1\big(\bar T; \F\big) \cong H^1\big(\bar S; \F\big)$.

We claim that, for small enough $\eps$, the map $f_\eps$ extends to a compactly supported $(C \cdot \eps)$-contracting map $F_{C\eps} \colon \bar T \to \SSS^{n-\ell}$ where $C > 0$ depends only on $S$ and $T$, but not on $\eps$.

The proof is by induction on the $k$-skeleta $T^{(k)} \subset T$ relative to $S$, where $k \in \N_0$.
If $S^{(0)} = T^{(0)}$, this is achieved by applying \cite[Lemma 3.2]{BH10} finitely many times.
If $S^{(0)} \subsetneq T^{(0)}$, the proof is a~bit more complicated as the cases $k=0$ and $k=1$ have to be treated simultaneously.
We refer the reader to the detailed argument in \cite[p.\ 474f.]{BH10}, which carries over to our case.

Let $c_1', \ldots c_\ell' \in H^1\big(\bar T; \F\big) \cong H^1\big(\bar S; \F\big)$ correspond to $c_1, \ldots, c_\ell \in H^1\big(\bar S; \F\big)$.
Using the naturality of the transfer map, \eqref{nonzero} implies
\begin{equation*}
 \big( (F_{C \eps})^*\big(\sigma^{n-\ell}\big) \cup (c'_1 \cup \cdots \cup c'_\ell) \big) \cap \pi_T^!(h_T) \neq 0 \in H_0\big(\bar T; \F\big).
\end{equation*}

Since $\eps > 0$ can be chosen arbitrarily small, this implies property \ref{shrinkingmap} of Definition \ref{largehom} for the class $h_T \in H_n(T; \F)$.
\end{proof}

\begin{Corollary} The complement $H_*(X ; \F) \setminus H_*^{\text{\rm tor-}\ell/\F\text{\rm -enl}}(X; \F) \subset H_*(X ; \F)$ is a linear subspace.
\end{Corollary}

\begin{proof} The most important part is to show that if $h_1 , h_2 \in H_*(X; \F)$ such that $h_1 + h_2 \in \smash{H_*^{\text{\rm tor-}\ell/\F\text{\rm -enl}}(X; \F)}$, then at least one of $h_1$ or $h_2$ is contained in \smash{$H_*^{\text{\rm tor-}\ell/\F\text{\rm -enl}}(X; \F)$}.
This follows because for the subcomplex $S \subset X$ witnessing \smash{$h_1 + h_2 \in H_*^{\text{\rm tor-}\ell/\F\text{\rm -enl}}(X; \F)$} we can assume by Lemma \ref{subcomplexinv} that both $h_1$ and $h_2$ lie in the image of $H_*(S; \F) \to H_*(X; \F)$.
The remaining details are left to the reader.
\end{proof}

\begin{Lemma} \label{invariance} Let $X$ and $X'$ be connected simplicial complexes with finitely presented fundamental groups.
Let $\eta \colon X \to X'$ be a map inducing an isomorphism of fundamental groups.
Then
\[
 (\eta_*)^{-1} \big( H_*^{\text{\rm tor-}\ell/\F\text{\rm -enl}}(X'; \F)\big) = H_*^{\text{\rm tor-}\ell/\F\text{\rm -enl}}(X; \F) .
\]
\end{Lemma}

\begin{proof}
We first show the inclusion ``$\subset$''.
This also implies the inverse inclusion, if $\eta$ is a~homotopy equivalence, by considering a~homotopy inverse of $\eta$.
Suppose $h \in H_*(X; \F)$ such that~\smash{$h' := \eta_*(h) \in H_*^{\text{\rm tor-}\ell/\F\text{\rm -enl}}(X'; \F)$}.
We need to show that \smash{$h \in H_*^{\text{\rm tor-}\ell/\F\text{\rm -enl}}(X; \F)$}.

There exists a finite connected subcomplex $S \subset X$ such that the inclusion $S \hookrightarrow X$ induces an isomorphism of fundamental groups and a class $h_S \in H_n(S; \F)$ which maps to $h \in H_n(X; \F)$.

Using Lemma \ref{subcomplexinv}, there exists a finite connected subcomplex $S' \subset X$ and $h'_{S'} \in H_n(S'; \F)$ satisfying the properties \ref{isofund} and \ref{shrinkingmap} of Definition \ref{largehom} for $h'$ such that, in addition, $\eta(S) \subset S'$ and~${\eta_*(h_S) = h'_{S'}}$.
Since $S$ is compact, the restriction $\eta|_{S} \colon S \to S'$ is, up to homotopy, $\lambda$-Lipschitz for some $\lambda > 0$.

Let $\eps > 0$, let $c'_1, \ldots, c'_\ell \in H^1\big(\bar S' ; \F\big)$, let $\pi_{S'} \colon \bar S' \to S'$ be a connected cover and let $f \colon \bar S' \to \SSS^{n-\ell}$ be a compactly supported $\tfrac{\eps}{\lambda}$-contracting map satisfying
\begin{equation} \label{forSprime}
 \big( f^*\big(\sigma^{n-\ell}\big) \cup (c'_1 \cup \cdots \cup c'_\ell) \big) \cap \pi_{S'}^!(h'_{S'}) \neq 0 \in H_{0}\big(\bar S' ; \F\big) .
\end{equation}
Since $\eta$ induces an isomorphism of fundamental groups, there is a unique connected cover $\pi_S \colon \bar S \to S$ such that $\eta$ lifts to a proper map $\bar \eta \colon \bar S \to \bar S'$ which induces an isomorphism of fundamental groups, hence an isomorphism $H^1\big(\bar S' ; \F\big) \cong H^1\big(\bar S; \F\big)$.
The lift $\bar \eta$ is still $\lambda$-Lipschitz continuous (recall that $\bar S$ is equipped with the canonical path metric).
Then $f \circ \bar \eta \colon \bar S \to \SSS^{n-\ell}$ is $\eps$-Lipschitz and compactly supported since $\bar \eta$ is proper.

Let $c_1 , \ldots, c_\ell \in H^1\big(\bar S; \F\big)$ correspond to $c'_1, \ldots, c'_\ell \!\in\! H^1\big(\bar S' ; \F\big)$.
By \eqref{forSprime} and since ${\eta_*(h_S)\! =\! h'_{S'}}$, we get
\begin{equation*}
 \bar \eta_* \big( \big( (f \circ \bar \eta)^*\big(\sigma^{n-\ell}\big) \cup (c_1 \cup \cdots \cup c_\ell) \big) \cap \pi_S^!(h_S) \big) \neq 0 \in H_0\big(\bar S' ; \F\big) .
 \end{equation*}
Since $\bar \eta_* \colon H_0\big( \bar S; \F\big) \to H_0\big(\bar S' ; \F\big)$ is an isomorphism, this shows that property \ref{shrinkingmap} of Definition~\ref{largehom} holds for $(S, h_S)$.
We conclude \smash{$h \in H_*^{\text{\rm tor-}\ell/\F\text{\rm -enl}}(X; \F)$}, which was our claim.

To show the reverse inclusion ``$\supset$'' in general, after replacing $X'$ by a homotopy equivalent simplicial complex, we may assume that $\eta \colon X \to X'$ is a simplicial inclusion.
Let \smash{$h \in H_*^{\text{\rm tor-}\ell/\F\text{\rm -enl}}(X; \F)$} and let $S \subset X$ and $h_S \in H_n(S; \F)$ be such that the properties \ref{isofund} and \ref{shrinkingmap} of Definition \ref{largehom} hold for $(S, h_S)$.
Then $S \hookrightarrow X'$ induces an isomorphism of fundamental groups and $H_n(S; \F) \to H_n(X'; \F)$ sends $h_S$ to $\eta_*(h)$.
This shows \smash{$\eta_*(h) \in H_*^{\text{\rm tor-}\ell/\F\text{\rm -enl}}(X'; \F)$}, finishing the proof of the reverse inclusion.
\end{proof}

For each finitely presented group $G$, we now define a subset
\[
 H_*^{\text{\rm tor-}\ell/\F\text{\rm -enl}}(G; \F) \subset H_*(G; \F)
\]
of {\em tor-$\ell / \F$-enleargeable classes} in the group homology of $G$ by setting
\[
 H_*^{\text{\rm tor-}\ell/\F\text{\rm -enl}}(G; \F) := H_*^{\text{\rm tor-}\ell/\F\text{\rm -enl}}(X; \F),
\]
where $X$ is an arbitrary connected simplicial complex representing $\sB G$.
By Lemma \ref{invariance}, this definition is independent of the choice of $X$.

\begin{Proposition} \label{homol_inv} A connected closed $\F$-oriented manifold $V$ is tor-$\ell/\F$ enlargeable if and only if the image of the fundamental class $[V]_{\F} \in H_n(V; \F)$ under the classifying map $V \to {\sf B}\pi_1(V)$ lies in \smash{$H_*^{\text{\rm tor-}\ell / \F\text{\rm -enl}}(\pi_1(V); \F)$}.
\end{Proposition}

\begin{proof} Choose a triangulation of $V$.
Then $V$ is tor-$\ell/\F$-enlargeable if and only if the fundamental class $[V]_\F \in H_n(V ; \F)$ is tor-$\ell/\F$-enlargeable.
This holds because the path metric induced by any Riemannian metric on $V$ is bi-Lipschitz equivalent to the simplicial path metric for the chosen triangulation on $V$ and each compactly supported $\eps$-Lipschitz map $\bar V \to \SSS^{n}$ of non-zero $\F$-toral degree is homotopic to a smooth compactly supported $2\eps$-Lipschitz map~${\bar V \to \SSS^{n-\ell}}$ of non-zero $\F$-toral degree.
Furthermore, the classifying map $V \to \sB \pi_1(V)$ induces an isomorphism of fundamental groups.
Thus the claim follows from Lemma \ref{invariance} for $X = V$ and~${X' = \sB \pi_1(V)}$.
\end{proof}

\begin{Example} \label{ex:torsion} Let $p \geq 2$ be a prime, let $\ell \geq 1$ and let $\tau_\ell \in H_\ell \big( (\Z/p)^\ell ; \F_p\big)$ be represented by the map
$
 \TTT^\ell = \sB \Z^{\ell} \to \sB (\Z/p)^{\ell}$
which is induced by the canonical projection $\Z^{\ell} \to (\Z/p)^{\ell}$.
Apparently, $\tau_\ell \in H_\ell\big( (\Z/p)^{\ell} ; \F_p\big)$ is $\F_p$-toral.

Next, let $(W_0, g)$ be a closed connected oriented Riemannian manifold of dimension $n$ where~$g$ has non-positive sectional curvature.
Let $\tilde W_0 \to W_0$ be the universal cover.
By the Cartan--Hadamard theorem, for any $x \in \tilde W_0$, the exponential map $T_x \tilde W_0 \to \tilde W_0$ is a diffeomorphism with $1$-Lipschitz inverse with respect to the lifted metric on $\tilde W_0$.
In particular, $W_0$ is a model for~$\sB \pi_1(W_0)$, and for all $\eps>0$ there exists a compactly supported $\eps$-contracting map $\tilde W_0 \to \SSS^n$ of degree $1$.
We claim that \smash{$\tau_\ell \times [W_0] \in H_{\ell + n}^{\text{\rm tor-}\ell/{\F_p}\text{\rm -enl}}\big((\Z/p)^{\ell} \times \pi_1(W_0) ; \F_p\big)$}.

To prove this, choose a simplicial complex $X$ representing $\sB (\Z/p)^{\ell}$.
Since $(\Z/p)^{\ell}$ is a finite group, one can choose $X$ finite in each dimension.
Let $K = X^{(2\ell)} \subset X$ be the $2\ell$-skeleton of $X$.
The inclusion $K \hookrightarrow X$ is $2 \ell $-connected.
In particular there is a unique $h_K \in H_\ell( K ; \F_p)$ mapping to $\tau_\ell$ under the map $H_*(K ; \F_p) \to H_*(X ; \F_p)$.

Choose some triangulation of $W_0$ and notice that the canonical path metric on $W_0$ is bi-Lipschitz equivalent to the metric induced by $g$.
Equip $X \times W_0$ with the product triangulation.
Then $K \times W_0$ with its product triangulation is a finite connected subcomplex of $X \times W_0$ and the class $h_K \times [W_0] \in H_{\ell + n}(K \times W_0 ; \F_p)$ maps to $\tau_\ell \times [W_0]$.
Furthermore, the inclusion $K \times W_0 \hookrightarrow X \times W_0$ induces an isomorphism of fundamental groups.

Let $\eps> 0 $, and let $f \colon \tilde W_0 \to \SSS^{n}$ be a compactly supported $\eps$-contracting map of degree $1$.
Then the composition
\[
 K \times \tilde W_0 \to \tilde W_0 \to \SSS^{n}
\]
is compactly supported, $\eps$-contracting (the first map is $1$-Lipschitz) and of non-zero $\F_p$-toral degree, as required.

As an application, assume that $\ell + n \geq 4$ and let an oriented manifold $V$ be obtained by surgery along $\ell$ pairwise disjoint embedded circles in $\TTT^{\ell} \times W_0$ representing the $p$-multiples of generators of the fundamental group $\pi_1\big({\TTT}^\ell\big) = \Z^\ell \subset \pi_1\big(\TTT^{\ell} \times W_0\big)$.

Then $\pi_1(V) \cong (\Z/p)^\ell \times \pi_1(W_0)$ and the classifying map $\varphi \colon V \to \sB \pi_1(V)$ is oriented bordant over $\sB \pi_1(V)$ to the map $\TTT^{\ell} \times W_0 = \sB \Z^{\ell} \times W_0 \to \sB (\Z/p)^{\ell} \times \sB \pi_1(W_0) = \sB \pi_1(V)$.
In particular, we obtain $\varphi_{*}([V]_{\F_p}) = \tau_\ell \times [W_0]$.
By Proposition \ref{homol_inv}, the manifold $V$ is tor-$\ell/\F_p$-enlargeable.

Note that $\varphi_*([V]) \in H_{\ell+n}(\pi_1(V); \Z)$ is a $p$-torsion class.
\end{Example}

\begin{Remark} One may define a quantitative version of tor-$\ell/\F$-enlargeability based on the notion of {\em $\tilde \square^{\perp}$-spread} of $\F$-toral homology classes rather than enlargeability, compare \cite[Section~7.1]{Gr21}.
\end{Remark}

\section[Atorality of PSC-3-manifolds]{Atorality of \yes-$\boldsymbol{3}$-manifolds} \label{nopsc}

We will show the following homological property of \yes-$3$-manifolds.

 \begin{Proposition} \label{prop:toralthree}
Let $\F$ be a field with $\charac( \F) \neq 2$.
Then closed orientable \yes-$3$-manifolds are $\F$-atoral.
\end{Proposition}

For the proof, we recall the following classification result of \yes-$3$-manifolds \cite{G-L1983}.

\begin{Proposition} \label{prop:decompthree} Closed connected orientable \yes-$3$-manifolds are diffeomorphic to connected sums of manifolds of the following type:
\begin{enumerate}[\myicon]\itemsep=0pt
\item closed connected $3$-manifolds with universal covers homotopy equivalent to $\SSS^3$,
\item copies of $\SSS^2\times \SSS^1$.
\end{enumerate}
In particular, they are not $\mathbb{Q}$-essential.
\end{Proposition}

Indeed, by the Kneser--Milnor prime decomposition \cite{Milnor1962}, every closed connected orientable $3$-manifold $V$ is diffeomorphic to the connected sum of manifolds of the form described in Proposition \ref{prop:decompthree} and aspherical manifolds, i.e., manifolds with contractible universal covers.

It remains to show that if $V$ is \yes, then aspherical summands cannot occur.
This non-trivial fact can be proved in several ways, see
\begin{enumerate}[\myicon]\itemsep=0pt
 \item \cite[Theorem 5.2]{S-Y1979}, based on the theory of {\em minimal hypersurfaces}, together with the resolution of Waldhausen's surface subgroup conjecture \cite{K-M2012};
\item \cite[Theorem 8.1]{G-L1983}, using the {\em index theory of the Dirac operator} on a (possibly non-compact) cover of $V$;
\item \cite[Theorem 3.4]{Ros1986}, showing that the {\em Rosenberg index} $\alpha(V) \in {\rm KO}_3(C^*_{\R} (\pi_1(V)))$ vanishes, together with the validity of the real Baum--Connes conjecture for $3$-manifold groups \cite{MOP08}, which in turn follows from the hyperbolization theorem of Thurston--Perelman;
\item \cite[Corollary 3.10.1\,(F$''$)]{Gr21}, estimating the {\em Uryson width} of the universal cover of $V$ and using the fact that groups which are quasi-isometric to trees are virtually free, see \cite[Theorem~7.19]{GdlH}.\footnote{For a remarkable generalization of this argument to manifolds in dimensions~4 and~5, see \cite{Ch-Li}.}
\end{enumerate}

\begin{Proposition} \label{ranks}
Let $V$ be a closed connected $n$-manifold, $n \geq 2$, with universal cover homotopy equivalent to $\SSS^n$.
Then $\dim_{\F} H^1 (V; \F) \leq 2$ for any field $\F$.
\end{Proposition}

\begin{proof} Let $G = \pi_1(V)$ and write $V = \tilde V / G$ where $\tilde V \to V$ is the universal cover.
Note that~$\tilde V$ is closed and orientable.
If the (free) action of $G$ on $\tilde V \simeq \SSS^n$ is orientation reversing, then $n$ is even and~${G = \Z/2}$ by the Lefschetz fixed point theorem.
In this case, we have $\dim_{\F} H^1(V ; \F) = 1$ if~${\charac(\F) = 2}$ and $\dim_{\F} H^1(V ; \F) = 0$ if $\charac(\F) \neq 2$.

So we can assume that the $G$-action on $\tilde V$ is orientation preserving.
A transfer argument shows that $H^*(V; \mathbb{Q}) \cong H^*(\SSS^n ; \mathbb{Q})$ and it remains to study the case when $\F$ is of finite characteristic~${p \geq 2}$.
Let $P < G$ be a Sylow $p$-subgroup and let $\pi \colon \tilde V / P \to \tilde V / G$ be the induced covering map.

Since $\tilde V$ is a homology $n$-sphere, it follows from \cite[Section XVI.9, Application 4]{C-E1956} that $P$ has periodic group cohomology of period $n+1$.
Using the classification of finite $p$-groups with periodic cohomology in \cite[Theorem XII.11.6]{C-E1956},
we conclude that $P$ is a cyclic $p$-group or generalized quaternion group.\footnote{Also compare \cite[Chapter IV.6]{Adem-Milgram}.}
Hence, $\dim_{\F} H^1\big(\tilde V / P; \F\big) = \dim_{\F} H^1(P ; \F) \leq 2$.
In fact, if~${p \geq 3}$, then $P$ is cyclic (since generalized quaternion groups are of even order), and hence~\smash{$\dim_{\F} H^1\big(\tilde V / P ; \F\big) \leq 1$}.

The induced map $\pi^* \colon H^*\big(\tilde V / G; \F\big) \to H^*\big(\tilde V / P ; \F\big)$ is injective since its composition with the cohomological transfer map \smash{$H^*\big(\tilde V / P; \F\big) \to H^*\big(\tilde V / G; \F\big)$} is multiplication by $\deg \pi = [G:P]$, which is prime to $p$.

From this follows the assertion of the proposition.
 \end{proof}

 \begin{Corollary} \label{cor:asphthree} Let $\charac(\F) \neq 2$ and let $V$ be a closed connected $\F$-toral $3$-manifold.
Then $V$ is diffeomorphic to a connected sum $K \sharp V'$ where $K$ is an aspherical $3$-manifold.
In particular, $V$ is $\mathbb{Q}$-essential.
\end{Corollary}

 \begin{proof}
By Proposition \ref{ranks}, and since for $\charac (\F) \neq 2$, every class in $H^1(V; \F)$ squares to $0$, the prime decomposition of $V$ must contain aspherical summands.
\end{proof}

The proof of Proposition \ref{prop:toralthree} follows from the combination of Proposition \ref{prop:decompthree} and Corollary~\ref{cor:asphthree}.

\begin{Remark} \label{rem:mixed}\quad
\begin{enumerate}[(a)]\itemsep=0pt
\item \label{leqzwei} It is easy to see that for $n \leq 2$ and all fields $\F$, all orientable connected $ \F$-toral $n$-manifolds are aspherical, hence $\mathbb{Q}$-essential, and are \no.
The $\F_2$-toral \yes-manifold $\mathbb{R}P^2$ shows that the orientability assumption cannot be dropped if $\charac(\F) = 2$.
\item \label{notinge4} Let $p$ be a prime, let $n \geq 4$ and let the oriented manifold $V^n$ be obtained by surgery along~$n$ pairwise disjoint embedded circles in $\TTT^n$ which represent the $p$-multiples of generators of the fundamental group $\pi_1({\TTT}^n) = \Z^n$.
Then $\pi_1(V) \cong (\Z/p)^n$, $V$ is $\F_p$-toral, but $V$ is not $\mathbb{Q}$-essential.
This shows that Corollary \ref{cor:asphthree} does not hold for $\F_p$-toral $n$-manifolds if $n \geq 4$.
\item The $\F_2$-toral \yes-manifold $\mathbb{R} P^3$ shows that Proposition \ref{prop:toralthree} does not hold for $p= 2$.
\end{enumerate}
\end{Remark}

\section[Torsion PSC-obstructions in higher dimensions]{Torsion \yes-obstructions in higher dimensions}

Using generalized minimal hypersurfaces ($\mu$-bubbles) and a Schoen--Yau-Schick descent, we propagate the \yes-obstruction from Section \ref{nopsc} to dimensions $\geq 4$.

We start with the following version of an overtorical band width inequality.
In the following, the scalar curvature of a Riemannian manifold $(V,g)$ is denoted by $\Sc_g(V) \colon V \to \R$.

\begin{Proposition} \label{prop:bandwidth} Let $\F$ be a field, let $2 \leq n \leq 7$ and let $0 \leq \ell < n$.
Also assume that one of the following conditions holds:
\begin{enumerate}[\myicon]\itemsep=0pt
 \item $\charac(\F) = 0$,
 \item $\charac(\F) = 2$ and $\ell \leq 2$,
 \item $\charac(\F) > 2$ and $\ell \leq 3$.
\end{enumerate}

Let $(B, g)$ be a compact oriented Riemannian $n$-manifold with boundary satisfying $\Sc_g(B) \geq n(n-1)$.
Let
$
 f \colon B \to [-1, +1] \times {\TTT}^{n- 1 -\ell}
$
be a smooth map satisfying $f(\partial B) \subset \{ -1, +1\} \times {\TTT}^{n- 1-\ell}$ such that there exist classes $c_1, \ldots, c_\ell \in H^{1}(B; \F)$ and a regular value $s \in [-1, +1] \times {\TTT}^{n- 1 -\ell}$ of~$f$ with\footnote{As $f(\partial B) \subset \{-1, +1\} \times {\TTT}^{n-1-\ell}$, the preimage $f^{-1}(s) \subset B$ is a closed smooth oriented submanifold, even if~${s \in \{ -1, +1\} \times \TTT^{n-1-\ell}}$.}
\[
\big\langle c_1 \cup \cdots \cup c_\ell, \big[f^{-1}(s) \big]_\F \big\rangle \neq 0 \in \F .
\]
Put $\partial_{\pm} B:= \partial B \cap f^{-1}\big( \{ \pm 1\} \times {\TTT}^{n- 1-\ell} \big) \subset \partial B$.
Then we have \smash{${\rm dist}_g(\partial_- B, \partial_+ B) \leq \frac{2\pi}{n}$}.
\end{Proposition}

\begin{Remark} The assumption $n \leq 7$ in Proposition \ref{prop:bandwidth} and Theorem \ref{thm:tornexist} can be weakened, see Remark \ref{rem:higherdim_mantra}.
\end{Remark}

\begin{proof}[Proof of Proposition \ref{prop:bandwidth}]
Without loss of generality, we can assume that $B$ is connected.
Suppose ${\rm dist}_g(\partial_- B , \partial_+B) > \frac{2\pi}{n}$.
We will show that this leads to a contradiction.

Up to homotopy relative to $\partial B$, $f$ respects some collar structures near $\partial_+ B \subset B$ and near $\{+1\} \times {\TTT}^{n-1-\ell} \subset [-1 , +1] \times {\TTT}^{n-1-\ell}$.
Let $s \in \{ +1\} \times {\TTT}^{n-1-\ell}$ be a regular value of $f|_{\partial_+ B} \colon \partial_+ B \to \{ +1\} \times {\TTT}^{n-1-\ell}$.
Then $s$ is also a regular value of $f$ and $f^{-1}(s) \subset \partial_+ B$.
So, by assumption and setting $\vartheta := c_1 \cup \cdots \cup c_{\ell} \in H^{\ell} (B; \F)$, we obtain
\begin{equation} \label{regvalbound}
 \big\langle \vartheta , \big[ ( f|_{\partial_+ B})^{-1} (s)\big]_\F \big\rangle \neq 0 \in \F.
\end{equation}
Let $\pi \colon [-1, +1] \times {\TTT}^{n-1-\ell} \to {\TTT}^{n-1-\ell}$ be the projection and let
$
 \omega^{n-1-\ell} \in H^{n-1-\ell}\big( {\TTT}^{n-1-\ell} ; \Z\big)
$
 be the cohomological fundamental class.
By \eqref{regvalbound} and Poincar\'e duality, we have
\begin{equation} \label{boundary}
 \big\langle \big( (\pi \circ f)^* \big(\omega^{n-1-\ell}\big) \big) \cup \vartheta , [ \partial_+ B]_\F \big\rangle \neq 0 \in \F.
\end{equation}

The triple $(B, \partial_- B, \partial_+ B)$ is a {\em band} in the sense of \cite[Definition 1.1]{Raede22}.
We claim that there exists a closed oriented smooth hypersurface $\Sigma^{n-1} \subset {\rm int}(B)$ which separates $\partial_- B$ and~$\partial_{+}B$ and admits a metric $g^{n-1}$ of positive scalar curvature.
Indeed, if this is not the case, then the $\frac{2\pi}{n}$-in\-equality \cite[Theorem 2.8]{Raede22}, see also \cite[Sections 3.6, 3.7 and 5.2]{Gr21}, implies that ${\rm dist}_g(\partial_- B , \partial_+B) \leq \frac{2\pi}{n}$, which contradicts our assumption.
Hence, there exists a hypersurface $\Sigma^{n-1}$ with the stated properties, and in particular it follows that $n \geq 3$.

Possibly after removing some components, the hypersurface $\Sigma^{n-1}$ is homologous to $\partial_+ B$ in~$B$, that is,
\begin{equation} \label{homologies}
 \big[\Sigma^{n-1}\big] = [\partial_+ B] \in H_{n-1}(B; \Z).
\end{equation}
Consider the composition
\[
 \Psi^{n-1} \colon\ \Sigma^{n-1} \stackrel{f|_{\Sigma^{n-1}}}{\longrightarrow} [ -1, +1] \times {\TTT}^{n-1-\ell} \stackrel{\pi}{\longrightarrow} {\TTT}^{n -1-\ell}.
\]
Furthermore, let $ \vartheta_{\Sigma^{n-1}} \in H^{\ell}(\Sigma^{n-1}; \F)$ be the restriction of $\vartheta \in H^{\ell}(B; \F)$.

By \eqref{boundary} and \eqref{homologies}, we have
\[
 \big\langle \big( \Psi^{n-1}\big)^* \big(\omega^{n-1-\ell}\big) \cup \vartheta_{\Sigma^{n-1}} , \big[\Sigma^{n-1}\big]_\F \big\rangle \neq 0 \in \F.
\]
Let $\ell' := \max\{ \ell, 2\}$.
Now we do a {\em Schoen--Yau--Schick descent} to construct closed oriented smooth Riemannian $i$-manifolds $\big(\Sigma^i, g^i\big)$, $\ell' \leq i \leq n-2$, satisfying $\Sc_{g^i} \big( \Sigma^i\big)> 0$ and
\[
 \Sigma^{\ell'} \subset \Sigma^{\ell'+1} \subset \dots \subset \Sigma^{n-2} \subset \Sigma^{n-1} ,
\]
together with maps $ \Psi^i \colon \Sigma^{i} \to {\TTT}^{i - \ell}$ so that, for $\ell' \leq i \leq n-2$,
\begin{enumerate}[(i)]\itemsep=0pt
 \item \label{dual} $\Sigma^{i} \subset \big(\Sigma^{i+1}, g^{i+1}\big)$ is a closed oriented stable minimal smooth hypersurface representing the dual of the cohomology class in $H^1\big(\Sigma^{i+1}; \Z\big)$ represented by the composition
 \[
 \Sigma^{i+1} \stackrel{\Psi^{i+1}}{\longrightarrow} {\TTT}^{i-\ell+1}\longrightarrow \SSS^1,
 \]
 where the second map is the projection onto the second factor in ${\TTT}^{i-\ell+1} = {\TTT}^{i-\ell} \times \SSS^1$.
 \item \label{restpsi} $\Psi^i$ is given by the composition
 \smash{$
 \Psi^{i} \colon\ \Sigma^{i} \stackrel{\Psi^{i+1}}{\longrightarrow} {\TTT}^{i-\ell+1} \longrightarrow {\TTT}^{i-\ell}
 $}
 where the second map is the projection onto the first factor in ${\TTT}^{i-\ell+1} = {\TTT}^{i-\ell} \times \SSS^1$.
 \item The metric $g^i$ is conformal to the metric on $\Sigma^i$ induced by $g^{i+1}$.
 \end{enumerate}
 For details of this construction, see \cite[Theorem 5.1]{S-Y1979} for $i =2$, the proof of \cite[Theorem 1]{SY} for~${i \geq 3}$, and the discussion in \cite{Sch98}.

 Let $\omega^{i-\ell} \in H^{i - \ell}\big( {\TTT}^{i-\ell}; \Z\big)$ be the cohomological fundamental class and let $\vartheta_{\Sigma^{i}}$ be the restriction of $\vartheta_{\Sigma^{i+1}}$ to $\Sigma^i$.
By induction, it follows from \ref{dual} and \ref{restpsi} that
\[
 \big\langle \big(\Psi^{\ell'}\big)^* \big(\omega^{\ell'-\ell}\big) \cup \vartheta_{\Sigma^{\ell'}} , \big[\Sigma^{\ell'}\big] \big\rangle \neq 0 \in \F.
 \]
 In particular, since $\big(\Psi^{\ell'}\big)^* \big(\omega^{\ell'-\ell}\big)$ and \smash{$\vartheta_{\Sigma^{\ell'}}$} are products of classes in \smash{$H^1\big( \Sigma^{\ell'} ; \F\big)$}, we see that \smash{$\big(\Sigma^{\ell'}, g^{\ell'}\big)$} is an $\ell'$-dimensional oriented smooth $\F$-toral \yes-manifold.
After we have passed to an appropriate connected component, we can assume that \smash{$\Sigma^{\ell'}$} is connected.

If $\ell' = 2$, this leads to a contradiction, see Remark \ref{rem:mixed}\,\ref{leqzwei}, so we can assume $\ell' = \ell \geq 3$.
If~${\charac(\F) = 0}$, this implies that $\Sigma^{\ell}$ is overtorical (see Remark \ref{rem:bemerkung}\,\ref{overtorical}), and we can do the Schoen--Yau--Schick descent again to show that $\Sigma^{\ell}$ is \no, contradiction.
If $\charac(\F)$ is finite, we have~${\ell = 3}$ and hence $\charac(\F) > 2$ by assumption.
So we get a contradiction to Proposition~\ref{prop:toralthree}.\looseness=1
\end{proof}

We obtain the following \yes-obstruction.

\begin{Theorem} \label{thm:tornexist} Let $\F$ be a field and $0 \leq \ell \leq n \leq 7$.
Let $V$ be a closed orientable tor-$\ell/\mathbb F$-enlargeable $n$-manifold and assume that one of the following conditions holds:
\begin{enumerate}[\myicon]\itemsep=0pt
 \item $\charac(\F) = 0$,
 \item $\charac(\F) = 2$ and $\ell \leq 2$,
 \item $\charac(\F) > 2$ and $\ell \leq 3$.
\end{enumerate}
Then $V$ is \no.
\end{Theorem}

\begin{proof} Suppose there is a counterexample $(V,g)$ with $\Sc_g(V) > 0$ of dimension $2 \leq n \leq 7$.
Choose an orientation of $V$.
By a scaling of $g$ we can assume $\Sc_g(V) \geq n(n-1)$.

If $\ell = n$, then a finite Riemannian covering $\bar V$ of $V$ is $\mathbb{F}$-toral.
If $\charac(\F) = 0$ then there exists a~map of non-zero degree $\bar V \to {\TTT}^n$, see Example \ref{rem:bemerkung}\,\ref{overtorical}.
Applying a Schoen--Yau--Schick descent as in the proof of Proposition \ref{prop:bandwidth}, this implies that $\bar V$ is \no.
If $\charac(\F) = 2$ and~${n = \ell = 2}$, then~$\bar V$ is \no~by Remark \ref{rem:mixed}\,\ref{leqzwei}.
If $\charac(\F) > 2$ and $n = \ell \leq 3$, then $\bar V$ is \no, by Proposition~\ref{prop:toralthree} for $n =3$, or Remark \ref{rem:mixed}\,\ref{leqzwei} for $n = 2$.
In summary, for $\ell = n$ we see that $\bar V$ and hence $V$ are~\no.

So we can assume that $\ell < n$.
Choose a smooth embedding
\[
 \psi \colon\ [-1,+1] \times \TTT^{n-1-\ell} \hookrightarrow \SSS^{n-\ell} \setminus \{\infty \},
\]
where $[-1,+1] \times \TTT^{n-1-\ell}$ carries the product of the standard metrics.
For $\delta \in (0,1]$, we obtain the smooth submanifold $A_\delta := \psi\big( [ - \delta , + \delta] \times \TTT^{n-1-\ell} \big) \subset \SSS^{n-\ell}$ with boundary
\[
 \partial A_\delta = \partial_- A_\delta \cup \partial_+ A_\delta , \qquad \partial_{\pm} A_\delta = \psi\big(\{ \pm \delta\} \times \TTT^{n-1-\ell}\big) .
\]
We equip $A_\delta$ with the Riemannian metric induced from $\SSS^{n-\ell}$.
Since $A_1 \subset \SSS^{n-\ell}$ is compact, the inverse $\psi^{-1} \colon A_1 \to [-1,+1] \times \TTT^{n-1-\ell}$ is $\Lambda$-Lipschitz for some $\Lambda > 0$.
Hence, for each $\delta \in (0,1]$, we have
\begin{equation} \label{estdist}
 2 \delta \leq \Lambda \cdot {\rm dist}_{A_\delta} ( \partial_- A_\delta, \partial_+ A_\delta) .
\end{equation}
By assumption, there exist a connected cover $\bar V \to V$, cohomology classes $c_1, \ldots , c_{\ell} \in H^1\big( \bar V; \F\big)$, an $\eps := \tfrac{n}{2 \pi \Lambda}$-contracting compactly supported map $f \colon \bar V \to \SSS^{n-\ell}$ and a regular value $s \in A_{\nicefrac{1}{2}} \subset \SSS^{n-\ell} \setminus \{ \infty\}$ of $f$ such that
\begin{equation} \label{regval}
 \big\langle c_1 \cup \cdots \cup c_{\ell}, \big[ f^{-1}(s) \big]_\F \big\rangle \neq 0 \in \F .
\end{equation}
Let $\delta \in \big(\tfrac{1}{2} ,1\big)$ be such that $f$ is transverse to $\partial A_\delta \subset \SSS^{n-\ell}$ and let
$
 B := f^{-1} ( A_\delta) \subset \bar V $.
This is a~smooth compact oriented submanifold of codimension $0$ with boundary
$
 \partial B = \partial_- B \cup \partial_+ B$, $ \partial_{\pm} B := f^{-1} (\partial_{\pm} A_\delta) $.
We equip $B$ with the Riemannian metric from $V$.
By our choice of $\varepsilon$ and using \eqref{estdist} together with $\delta > \tfrac{1}{2}$, we have
\begin{equation} \label{breite}
 {\rm dist}_B( \partial_- B, \partial_+ B) \geq \frac{2 \pi \Lambda}{n} \cdot {\rm dist}_{A_\delta} (\partial_- A_\delta, \partial_+ A_\delta)
 \geq \frac{2 \pi \Lambda}{n} \cdot \Lambda^{-1} \cdot 2\delta > \frac{2\pi}{n} .
\end{equation}
Let $c'_1, \ldots , c'_\ell \in H^1(B; \F)$ denote the restrictions of $c_1, \ldots , c_\ell \in H^1\big(\bar V; \F\big)$.
We obtain by \eqref{regval}
\[
 \big\langle c'_1 \cup \cdots \cup c'_{\ell}, \big[ f^{-1}(s) \big]_{\F} \big\rangle \neq 0 \in \F.
\]
Using $\Sc_g(B) \geq n(n-1)$ and \eqref{breite}, this contradicts Proposition \ref{prop:bandwidth} since $A_{\delta}$ is diffeomorphic to~${[-1, +1] \times \TTT^{n-1-\ell}}$.
\end{proof}

\begin{Remark} \label{rem:higherdim_mantra}\qquad
 \begin{enumerate}[(a)]\itemsep=0pt
\item \label{mantra1} If $n = 8$, then \cite{Smale} implies that the $\mu$-bubble construction underlying the proof of \cite[Theorem 2.8]{Raede22} can be performed after a small $C^{\infty}$-perturbation of $g$.
In view of the recent results in \cite{CMS}, this can also be achieved for $n = 9, 10$.
Higher dimensions $n \geq 11$ may possibly be covered by the results in \cite{Loh2018}.
Thus the condition $n \leq 7$ in Theorem \ref{thm:tornexist} can be weakened accordingly.
\item The real projective spaces $\R P^n$ are $\F_2$-toral for all $n$, are \yes~for $n \geq 2$, and orientable for all odd $n \geq 3$.
Thus, if $\charac(\F) = 2$, then neither the orientability of $V$ nor the assumption~${\ell \leq 2}$ in Theorem \ref{thm:tornexist} can be dropped.
\end{enumerate}
\end{Remark}

\begin{Example} \label{ex:counterexample} Let $4 \leq n \leq 7$, let $0 \leq \ell \leq n $ and let the orientable manifold $W$ be obtained from~${{\TTT}^n = {\TTT}^{\ell} \times {\TTT}^{n-\ell}}$ by killing the $2$-multiples of generators of $\pi_1\big({\TTT}^\ell\big) = \Z^\ell$ by surgery along pairwise disjoint embedded circles in ${\TTT}^n$.
Let $V$ be the connected sum of $W$ and $\mathbb{C} P^{2} \times \SSS^{n-4}$.

Then $V$ is an orientable tor-$\ell / \F_{\!2}$-enlargeable manifold (compare Example \ref{ex:torsion}) which is~\yes, if and only if $\ell \geq 3$.

The \no~property for $\ell \leq 2$ follows from Theorem \ref{thm:tornexist} while the \yes~property for $3 \leq \ell \leq n$ is shown as follows:
We have $\pi_1(V) = (\Z/2)^{\ell} \times \Z^{n-\ell}$ and the classifying map $\varphi \colon V \to \sB \pi_1(V)$ sends the integral fundamental class of $V$ to the cross product of the toral homology classes in~${H_{\ell}\big((\Z/2)^{\ell} ; \Z\big)}$ and $H_{n-\ell}\big(\Z^{n-\ell}; \Z\big)$.
Since $\ell \geq 3$, by \cite[Theorem A]{Jo2004} the toral homology class in~${H_{\ell}\big((\Z/2)^{\ell} ; \Z\big)}$ can be represented by a closed oriented \yes~$\ell$-manifold, and so the same holds for $\varphi_*([V]) \in H_n(\pi_1(V) ; \Z)$.
Using the fact that the universal cover of $V$ is non-spin, Theorem~\ref{prop:nonspin} implies that $V$ is \yes.
\end{Example}

The following fundamental conjecture remains open.

\begin{Conjecture} \label{main_conj}
Let $\F$ be a field of odd characteristic, let $n \geq 4$ and let $4 \leq \ell \leq n$.
Then each closed tor-$\ell/\mathbb F$-enlargeable $n$-manifold is \no.
\end{Conjecture}

Recall that Corollary \ref{cor:asphthree} is specific to $3$-dimensional manifolds, see Remark \ref{rem:mixed}\,\ref{notinge4}, and thus our proof of Theorem \ref{thm:tornexist} does not generalize to $\ell \geq 4$ if $\charac(\F) > 2$.

Table \ref{tab:conditions} collects our knowledge about the existence of closed orientable tor-$\ell / \F$-enlargeable \yes-manifolds of dimension $n \leq 7$ (or more general $n$ as in Remark \ref{rem:higherdim_mantra}\,\ref{mantra1}).

\begin{table}[h]\renewcommand{\arraystretch}{1.2}
\centering
\begin{tabular}{|c|c|c|c|}
\hline
 & $\charac(\F)$ = 0 & $\charac(\F) = 2$ & $\charac(\F) > 2$ \\
\hline
$\ell \leq 2$ & no & no & no \\
\hline
$\ell = 3$ & no & yes & no \\
\hline
$\ell > 3$ & no & yes & \textcolor{red}{open} \\
\hline
\end{tabular}

\caption{\label{tab:conditions}Existence of \yes-manifolds.}
\end{table}

\section[PSC-products]{\yes-products of \no-manifolds}

Given a closed connected oriented manifold $V$, we denote by
\[
 \Phi_{V} = \varphi_*([V]) \in \Ho_*( \pi_1(V);\Z)
\]
 the image of the homological fundamental class of $V$ under the classifying map $\varphi \colon V \to \sB \pi_1(V)$.

As an application of our discussion we obtain many \no-manifolds whose Cartesian products are \yes.

\begin{Theorem}\label{thm:plenty}
Let $V_1$, $V_2$ be closed oriented connected manifolds such that
\begin{enumerate}[\myicon]\itemsep=0pt
 \item $\dim (V_1 \times V_2) \geq 5$;
 \item the universal cover of $V_1 \times V_2$ is non-spin;
 \item $\Phi_{V_1} \in \Ho_*(\pi_1(V_1) ; \Z)$ and $\Phi_{V_2} \in H_*(\pi_1(V_2); \Z)$ are of finite, coprime order.
\end{enumerate}
Then $V_1 \times V_2$ is \yes.
\end{Theorem}

\begin{proof} The class $\Phi_{V_1 \times V_2} \in \Ho_*(\sB ( \pi_1(V_1) \times \pi_1(V_2)); \Z)$ is the image of $\Phi_{V_1} \otimes \Phi_{V_2}$ under the cross-product homomorphism
\[
 \Ho_*( \pi_1(V_1); \Z) \otimes \Ho_* ( \pi_1(V_2); \Z) \to \Ho_* ( \pi_1(V_1) \times \pi_1(V_2); \Z) .
\]
Hence, this class is equal to $0$, since $\Phi_{V_1}$ and $\Phi_{V_2}$ are elements of coprime order in the left-hand tensor factors.
Now the claim follows from Theorem \ref{prop:nonspin}.
\end{proof}

\begin{Example} \label{interesting}
Let $p_1$, $p_2$ be different odd primes, let $4 \leq n \leq 7$, let $1 \leq \ell \leq 3$ and let~$V_1$,~$V_2 $ be tor-$\ell /\F_{p_i}$-enlargeable manifolds of dimension $n$ obtained from ${\TTT}^{n} = {\TTT}^{\ell} \times {\TTT}^{n-\ell}$ by killing the $p_i$-multiples of the $\ell$ generators in $\pi_1\big({\TTT}^\ell\big) =\Z^{\ell}$ by surgery, see Example \ref{ex:torsion}.
Replace $V_1$ by~${V_1 \sharp \big( \mathbb{C}P^2 \times \SSS^{n-4} \big)}$ in order to make sure that the universal cover of $V_1 \times V_2$ is non-spin.
Then
\begin{enumerate}[\myicon]\itemsep=0pt
 \item $V_1$ and $V_2$ are \no~by Theorem \ref{thm:tornexist}.
 \item $V_1 \times V_2$ is \yes~by Theorem \ref{thm:plenty} and the final statement in Example \ref{ex:torsion}.
\end{enumerate}
\end{Example}

\begin{Remark}\quad
\begin{enumerate}[(a)]\itemsep=0pt
\item Examples \ref{interesting} elaborates on the $\mathcal{SYS} \times \mathcal{SYS}$-remark from \cite{Gr2018}.
They are the first examples of their kind in dimensions greater than $4$ and with non-spin universal covers.
\item Complex surfaces of odd degree $\geq 5$ in $\mathbb{C}P^3$ are \no~by Seiberg--Witten theory, while their products are \yes, as are all simply connected non-spin $n$-manifolds for $n\geq 5$, see~\cite{GL}.
\item A referee pointed out the following class of spin examples.
Recall the graded coefficient ring of connective real $K$-homology,
\[
 {\rm ko}_* \cong \Z[ \eta, \omega, \mu] / \big( 2 \eta, \eta^3, \omega \eta, \omega^2 - 4 \mu\big), \qquad |\eta| = 1,\qquad |\omega| = 4,\qquad |\mu|=8.
\]
By \cite[Theorem A]{Stolz92}, a simply connected closed spin manifold of dimension $n \geq 5$ is~\yes~if and only if the Hitchin index \cite{Hitchin} vanishes, $\alpha(M) = 0 \in {\rm ko}_n$.
Since ${\rm ko}_*$ has zero divisors, one obtains simply connected closed spin \no~manifolds $V_1$ and $V_2$ whose products are~\yes.

For example, it is well known, see \cite[p.\ 44]{Hitchin}, that for $n \equiv 1, 2 \mod(8)$, there are exotic spheres $\Sigma^n$ with $\alpha(\Sigma^n) \neq 0 \in {\rm ko}_n \cong \Z/2$.
Let $K^4$ be the K3 surface.
It is simply connected and satisfies $\alpha(K^4) = \omega \in {\rm ko}_4 \cong \Z$.
Since $\eta \cdot \omega = 0$ and the Hitchin index is multiplicative, the products $\Sigma^n \times K^4$ are \yes, whereas $\Sigma^n$ and $K^4$ are \no.
\end{enumerate}
\end{Remark}

\section{Four conjectures} \label{fourconjectures}

\begin{Conjecture}[short/long neck torality inequality\footnote{This is motivated by Cecchini's long neck inequality \cite{Ce2020}.}]
\label{shortneck}
Let $0 \leq \ell < n$ be natural numbers, let~$(V,g)$ be a compact oriented Riemannian $n$-manifold with boundary and let
$f \colon  V\to \SSS^{n-\ell}$
be a smooth {\it area decreasing} map\footnote{That is ${\rm area}(f(\Sigma)) \leq {\rm area}(\Sigma)$ for all smooth surfaces $\Sigma\subset X$.} which is constant on the boundary $\partial V\subset V $, say $f(\partial V)=\{ \infty \}\subset \SSS^{n-\ell}$.
Suppose $\charac(\F) \neq 2$ and there are classes $c_1, \ldots, c_\ell \in H^1(V; \F)$ such that for some regular value $s \in \SSS^{n-\ell} \setminus \{ \infty\}$ of $f$ we get
\[
 \bigl\langle c_1 \cup \cdots \cup c_\ell, \bigl[ f^{-1}(s)\bigr] \bigr\rangle \neq 0 \in \F .
\]
Setting $k := n - \ell$, suppose that for all $x \in \operatorname{supp}( {\rm d}f) = \{ x \in V \mid {\rm d}f(x) \neq 0\}$, we have
\[
 \Sc_g(V, x) \geq k(k-1) .
\]
If the scalar curvature of $V$ is bounded from below by $\sigma > 0$, i.e.,
\[
 \Sc_g(V,x)\geq \sigma>0 \qquad \textrm{for } x\in V,
 \]
then the distance from the support of ${\rm d}f$ to the boundary of $V$ satisfies
\[
 {\rm dist}_g ( \operatorname{supp}({\rm d}f), \partial V) < {\rm const} \cdot \sqrt{\tfrac {1}{\sigma}},
 \]
 possibly for \smash{${\rm const}=\pi \sqrt{\frac{n-1}{n}}$} where this constant is optimal.
 \end{Conjecture}

 \begin{Remark} There is a significant mismatch between distance and area domination with positive scalar curvature.
Namely, given an orientable $n$-manifold $V=(V,g)$ with $\Sc_g(V) \geq \sigma>0$, then, for all $\varepsilon >0$, there exists a metric $g_{\varepsilon}$ on $V$, such that
\begin{enumerate}[\myicon]\itemsep=0pt
\item $ V_{\varepsilon}=(V,g_{\varepsilon})$ admits an {\it area decreasing} diffeomorphism $V\to V_{\varepsilon }$;
\item no Riemannian manifold $(W,h)$ which $1$-Lipschitz dominates $V_\eps$, i.e., which admits a proper $1$-Lipschitz map of non-zero degree $W \to V_\eps$, has $\Sc_h \geq \varepsilon$.
\end{enumerate}
The simplest example is where $V=\SSS^2\times \TTT^{n-2}$ and $V_\varepsilon =\SSS^2_{\varepsilon, D}\times \TTT_\varepsilon^{n-2}$, where $\SSS^2_{\varepsilon, D}$ is the (smoothed) boundary of the cylinder with base of radius $\varepsilon$ and height \smash{$D=10\sqrt{ \tfrac {1}{\varepsilon}}$}, i.e.,
\[
 \SSS^2_{\varepsilon, D}=\partial \big( B^2(\varepsilon)\times [0,D] \big),
 \]
$\TTT^{n-2} = \SSS^1 \times \cdots \times \SSS^1$ is the $(n-2)$-torus with the standard metric and $\TTT_\varepsilon^{n-2}$ is the product of circles with radius $\varepsilon$.
Here, the existence of an area-decreasing diffeomorphism $V\to V_{\varepsilon, }$ (for~${\varepsilon \ll 0.01}$) is elementary, while the non-existence of $1$-Lipschitz domination with $Sc \geq \varepsilon$ follows from the~$\frac {2\pi}{n}$-inequality as $10 > 2\pi$, see \cite[Section 3.6]{Gr21} and \cite[Theorem 2.8 and Proposition~2.21]{Raede22}.\footnote{The case of $n \leq 7$ is classical, while the case $8 \leq n \leq 10$ uses the methods of \cite{Smale} and \cite{CMS}, respectively.
For $n \geq 11$, see \cite[Theorem 4.6]{S-Y2017} and \cite{Loh2018}.
If $V_\eps$ is spin, then the full $\frac {2\pi}{n}$-inequality is proved in \cite[Theorem 1.4]{Ze2019}, \cite[Corollary 1.5]{Ze2020}, \cite[Theorem D]{Ce2020} and \cite[Theorem 1.3]{G-X-Y2020}.}

Now, for an arbitrary $(V,g)$, let $\TTT^{n-1}\subset V $ be a smoothly embedded torus and let us modify the metric $g$ in a small neighborhood $\TTT^{n-1}\times [-\delta, +\delta]\subset V$ of the torus by $\varepsilon'$-shrinking $g$ along the torus and \smash{$\sqrt \frac {1}{\varepsilon'}$}-stretching $g$ normally to the torus.
For sufficiently small $\eps' = \eps'(g, \eps) > 0$, this decreases the areas of the surfaces in $V$, while the $\frac {2\pi}{n}$-inequality prevents $1$-Lipschitz domination of $V_\varepsilon$ with $\Sc \geq \varepsilon$.
\end{Remark}

\begin{Conjecture}[tor-non-enlargeability of foliations with $\Sc>0$] \label{foliations}
 Let $V$ be a closed $n$-manifold and $\mathscr F$ be a smooth foliation on $V$ equipped with a leafwise Riemannian metric $g_{\mathscr F}$, i.e., $g_{\mathscr F}$ is a smooth bundle metric on the integrable subbundle $T \mathscr F \subset TV$ tangent to the leaves of $\mathscr F$.

Let $\Sc(\mathscr F) \colon V \to \R$ be defined as $\Sc(\mathscr F, x) = \Sc_{g_{\mathscr F}}({\mathscr F}_x , x)$ where ${\mathscr F}_x$ is the leaf passing through $x \in V$ equipped with the Riemannian metric induced by $g_{\mathscr F}$.

If $V$ is tor-$\ell / \mathbb F$-enlargeable for some non-negative integer $\ell$ and some field $\F$ of $\charac(\F) \neq 2$, then $\Sc(\mathscr F) \ngtr 0$.
 \end{Conjecture}

For foliations of codimension $1$, we have the following partial answer.

\begin{Proposition} Let $\F$ be a field, let $1 \leq n \leq 6$ and let $0 \leq \ell \leq n$.
Also assume that one of the following conditions holds:
\begin{enumerate}[\myicon]\itemsep=0pt
 \item $\charac(\F) = 0$,
 \item $\charac(\F) = 2$ and $\ell \leq 2$,
 \item $\charac(\F) > 2$ and $\ell \leq 3$.
\end{enumerate}
Let $V$ be a closed orientable tor-$\ell/\F$-overtorical\footnote{See Example \ref{ex:toral}\,\ref{ex:overtorical}.} $n$-manifold.

Then $V$ does not admit a smooth orientable codimension-$1$ foliation $\mathscr F=(\mathscr F, g_{\mathscr F})$ together with a leafwise Riemannian metric such that $ \Sc (\mathscr F)>0$.
\end{Proposition}

\begin{proof}[Sketch of proof] Suppose $V$ carries a smooth leafwise Riemannian metric $g_{\!\mathscr F}$ with ${\Sc_{g_{\!\mathscr F}} (\mathscr F)\!\! > \!\!0}$.
After scaling $g_{\mathscr F}$, we can assume that $\Sc_{g_{\mathscr F}} (\mathscr F) > (n+1)n$.
Let $V^* \to V $ be the (one-dimensional for $\mathrm{codim} (\mathscr F)=1$) {\it Connes bundle} of $\mathscr F$, which in this case is equivalent to the trivial $\R$-bundle~${V\times \mathbb R \to V}$, see \cite[Sections $1\frac{7}{8}$ and $9\frac {2}{3}$]{Gr1996} and \cite[Section 6.5.2]{Gr21}.
Recall \cite[p.\ 9]{Gr1996} that the total space $V^*$ of this bundle has an induced codimension-$2$ foliation $\mathscr F^*$ whose leaves correspond to leaves in $\mathscr F$, and that $V^*$ has a Riemannian metric $\tilde g$ such that
\begin{enumerate}[\myicon]\itemsep=0pt
 \item on each leaf of $\mathscr F^*$, the metric induced by $\tilde g$ coincides with the metric induced by $g_{\mathscr F}$ on the corresponding leaf of $\mathscr F$,
 \item the restriction of $\tilde g$ to $(T \mathscr F^*)^{\perp} \subset TV^*$ is {\em transversally Riemannian}, i.e., it is a bundle metric which is parallel with respect to the Bott connection on $(T \mathscr F^*)^{\perp} \to V^*$ along $\mathscr F^*$.
\end{enumerate}

For $\lambda > 0$, we define a new Riemannian metric $\tilde g_{\eps}$ on $V^* \cong V \times \R$ as the bundle metric on
\[
{TV^* = T\mathscr F^* \oplus (T\mathscr F^*)^{\perp}}
\]
 equal to
\[
 \tilde g_\lambda = \tilde g|_{T\mathscr F^*} + \lambda^2 \tilde g|_{(T\mathscr F^*)^{\perp}}.
\]
The discussion in \cite[p.\ 315]{Gr21} shows that there exist $\lambda \gg 0$ such that
\begin{enumerate}[\myicon]\itemsep=0pt
\item $\Sc_{\tilde g_{\lambda}} (V \times [-1, +1]) \geq (n+1)n$,
\item ${\rm dist}_{\tilde g_\lambda}( V \times \{-1\}, V\times \{+1\}) > \frac{2\pi}{n+1}$.
\end{enumerate}
Since $V$ is tor-$\ell/\F$-overtorical, this contradicts Proposition \ref{prop:bandwidth}.
\end{proof}

The above suggests the following.

\begin{Conjecture}[domination of foliations] If a closed orientable manifold $V$ admits a smooth foliation $\mathscr F$ with $\Sc(\mathscr F) > 0$, then $V$ can be dominated by an orientable manifold $W$ with $\Sc(W)>0$.\footnote{This is hardly to be expected from {\it leaf-wise} smooth Riemannian foliations $\mathscr F$ which are (only) {\it continuous} transversally to the leaves.
However, the standard examples of such $\mathscr F$ -- the stable/unstable foliations of {\it Anosov systems} -- are easily shown to be \no.}
\end{Conjecture}

We conclude by bringing together Conjectures \ref{shortneck} and \ref{foliations} as follows, compare \cite{S-W-Z2021}.

\begin{Conjecture}[foliated short/long neck torality inequality]
Let $V$ be a compact oriented Riemannian $n$-manifold with boundary and $\mathscr F$ be a smooth $m$-dimensional foliation on $V$.

Let $\charac(\F) \neq 2$ and let $f \colon V \to \SSS^{n-\ell}$ be a smooth map which {\it decreases the areas of smooth surfaces in the leaves of} $\mathscr F$ such that $f(\partial V)=\{ \infty\} \subset \SSS^{n-\ell}$ and such that there are classes~${c_1 , \ldots, c_\ell \in H^1(V ; \F)}$ with
\[
 \big\langle c_1 \cup \cdots \cup c_{\ell}, \big[f^{-1}(s)\big] \big\rangle \neq 0 \in \F
\]
for some regular value $s \in \SSS^{n-\ell} \setminus \{ \infty\}$.

Setting $k := n-\ell$, let
\[
 \Sc(\mathscr F,x)\geq\min \left\{ k(k-1), m(m-1)\right\}
\]
for all $x\in V$ with $x \in \operatorname{supp}( {\rm d}f_{\mathscr F}) = \{ x \in V \mid {\rm d}f(x)|_{T \mathscr F} \neq 0\}$.

If the scalar curvature of $\mathscr F$ is bounded from below by $\sigma > 0$, i.e.,
\[
 \Sc(\mathscr F,x)\geq \sigma>0 \qquad\text{for}\quad x\in V,
 \]
then {\it the leaf-wise distance} from the support of ${\rm d}f$ on $T\mathscr F$ to the boundary of $V$ satisfies
\[
 {\rm dist}_{\mathscr F}(\operatorname{supp}({\rm d}f_{\mathscr F}), \partial V)=\inf_{x\in \operatorname{supp} {({\rm d}f_{\mathscr F})}} {\rm dist}_{ \mathscr F_x} (x, \mathscr F _x\cap \partial X) < {\rm const} \sqrt{ \tfrac{1}{\sigma}},
\]
where $\mathscr F _x\subset V$ is the leaf of $\mathscr F$ through $x$.
\end{Conjecture}

\appendix

 \section[PSC-metrics on manifolds with non-spin universal covers]{\yes-manifolds with non-spin universal covers}

In \cite{GL}, it was proved that closed simply connected non-spin manifolds of dimension $\geq 5$ are~\yes.
Here we prove the corresponding result for non-simply connected manifolds.

\begin{Theorem} \label{prop:nonspin} Let $n \geq 5$ and let $V$ be a closed connected oriented smooth $n$-manifold with non-spin universal cover.
Let $\varphi \colon V \to {\sf B} \pi_1(V)$ be the classifying map and assume that there exists a map $f \colon V' \to {\sf B}\pi_1(V)$ where $V'$ is a closed oriented $n$-dimensional \yes-manifold such that
\[
 \varphi_*([V]) = f_{*} ( [ V' ]) \in H_n(\pi_1(V) ; \Z) .
\]
$($This includes the case $\varphi_*([V]) = 0$.$)$

Then $V$ is \yes.
\end{Theorem}

The proof is based on the following generic construction of \yes-manifolds, compare \cite{Stolz92} and \cite{Fuehring}.
Consider the action of $G := \SU(3) \rtimes \Z/2$ on $\C P^2$ where $\SU(3)$ acts on $\C P^2$ by matrix multiplication and $\Z/2$ acts by complex conjugation.
This action is transitive, its stabilizer group at $[0:0:1] \in \C P^2$ is equal to $H = {\rm S}({\rm U}(2) {\rm U}(1)) \rtimes \Z/2 < G$, and it is isometric with respect to the Fubini study metric on $\C P^2$.
Hence, we obtain a fibre bundle
\begin{equation} \label{universal}
 \C P^2 = G/ H \hookrightarrow \sB H \to \sB G
\end{equation}
with structure group $G$.
Let $B$ be a~closed oriented smooth manifold, let $X$ be a CW complex and let $
 f \times \omega \colon\ B \to X \times \sB G
$
be a continuous map.
The map $\omega \colon B \to \sB G$ classifies a smooth fiber bundle $\C P^2 \hookrightarrow M_{\omega} \stackrel{\pi}{\longrightarrow} B$ with structure group $G$, which is unique up to fibre bundle diffeomorphism.
The total space $M_{\omega}$ is a closed oriented smooth manifold and it is \yes, by the O'Neill formula.
Furthermore, the assignment $(f \times \omega \colon B \to X \times BG) \mapsto (f \circ \pi \colon M_\omega \to X)$ respects the oriented bordism relation, and we obtain a group homomorphism which is natural in $X$,
\begin{equation} \label{psi}
 \psi(X)_* \colon\ \Omega^{\SO}_{*-4} (X \times \sB G) \to \Omega^{\SO}_{*}(X).
\end{equation}
Let $\mathscr{T}_*(X) := \im ( \psi(X)_*) \subset \Omega^{\SO}_{*}(X)$ and let $ \Omega^{\SO, +}_*( X) \subset \Omega_*^{\SO}(X)$ be the subgroup consisting of oriented bordism classes $[M \to X]$ where $M$ is \yes.
We obtain
$
 \mathscr{T}_*(X) \subset \Omega^{\SO,+}_*(X)$.
Now consider the {\em homological orientation}, which is natural in $X$,
\begin{equation} \label{homor}
 u(X)_* \colon\ \Omega^{\SO}_*(X) \to H_*(X; \Z), \qquad [f \colon M \to X] \mapsto f_*([M]) \in H_*(X; \Z).
\end{equation}
For dimension reasons, $\pi \colon M_\omega \to B$ sends the fundamental class of $M_\omega$ to zero, hence
$
 \mathscr{T}_*(X) \subset \ker u(X)_*$.
We will next show that $2$-locally, the reverse inclusion holds as well.

For an abelian group $A$, we will use the shorthand $A_{(2)} = A \otimes \Z_{(2)}$ for the $2$-localisation.

 \begin{Proposition} \label{crucial} We have
$
 \mathscr{T}_*(X)_{(2)} = (\ker u(X)_*)_{(2)}$.
 In particular, for all $x \in \ker u(X)_*$ there is an odd $k \in \N$ with $kx \in \mathscr{T}_*(X)$.
\end{Proposition}

\begin{proof}
Only the inclusion ``$\supset$'' remains to be shown.
It is well known, see \cite[Chapter 2]{CF}, that the Atiyah--Hirzebruch spectral sequence for $2$-local oriented bordism
\[
 E^2_{p,q} \cong H_p\big(X; \Omega^{\SO}_q\big)_{(2)} \Longrightarrow \Omega^{\SO}_{n}(X)_{(2)}
 \]
collapses at $E^2$.
There is a decreasing filtration
\[
 \Omega^{\SO}_{n}(X)_{(2)} = \mathscr{F}_{n,0} \supset \cdots \supset \mathscr{F}_{0,n} = H_0(X ; \Omega_n)_{(2)}
\]
with \smash{$E^{2}_{p,q} \cong \mathscr{F}_{p,q} / \mathscr{F}_{p-1,q+1}$}.
Furthermore, the usual $\Omega^{\SO}_*$-module structure on $\Omega^{\SO}_*(X)$ induces maps
$
 \Omega^{\SO}_{\mu} \times \mathscr{F}_{p, q} \mapsto \mathscr{F}_{p, q+\mu}
$
so that the isomorphism $E^2_{p,*} \cong H_p\big(X; \Omega^{\SO}_*\big)$ is $\Omega^{\SO}_*$-linear.

The homological orientation $u(X)_*$ corresponds $2$-locally to the projection
\[
 \Omega^{\SO}_{n}(X)_{(2)} = \mathscr{F}_{n,0} \to \mathscr{F}_{n,0} / \mathscr{F}_{n-1,1} = H_n\big(X; \Omega^{\SO} _{0}\big)_{(2)} = H_n(X; \Z_{(2)}).
\]
Thus, by the collapsing of the Atiyah--Hirzebruch spectral sequence, $u(X)_*$ is $2$-locally surjective.

Let $\Tor_* \subset \Omega^{\SO}_*$ be the torsion subgroup.
The short exact sequence
\[
 0 \to \Tor_* \to \Omega^{\SO}_* \to \Omega^{\SO}_* / \Tor_* \to 0
\]
splits, and hence it induces short exact sequences
\begin{equation} \label{ses}
 0 \to H_p(X; \Tor_q)_{(2)} \hookrightarrow E^2_{p,q} \longrightarrow H_p\big(X; \Omega^{\SO}_q / \Tor_q\big)_{(2)} \to 0 .
\end{equation}
It is well known that $\Omega^{\SO}_*$ contains no odd torsion, so we have $H_p(X; \Tor_q)_{(2)} = H_p(X; \Tor_q)$.

The proof of Proposition \ref{crucial} is complete once we show that every element in $E^2_{p,q}$ with $p + q = n$ and $q > 0$ can be represented by an element in $\mathscr{T}_n(X)_{(2)} \cap \mathscr{F}_{p,q}$.
According to \eqref{ses}, we decompose this claim into the following two assertions.

\begin{Assertion}\label{Assertion1} Let $q > 0$.
Then the image of $\mathscr{T}_n(X)_{(2)} \cap \mathscr{F}_{p,q}$ in $E^2_{p,q}$ is mapped surjectively to $H_p\big(X; \Omega_q^{\SO} / \Tor_q \big)_{(2)}$.
\end{Assertion}

To show this, consider a simple tensor
\[
 h \otimes \xi \in H_p\big(X; \Omega_q^{\SO} / \Tor_q \big)_{(2)} = H_p(X; \Z_{(2)}) \otimes \Omega_q^{\SO} / \Tor_q .
\]
There exists an odd $k \in \N$ such that $k\cdot h \in H_p(X; \Z_{(2)})$ lifts to a bordism class $ [ f \colon M^p \to X] \in \Omega_p^{\SO}(X)$.
Furthermore, by constructing appropriate $\C P^2$-bundles, see \cite[Proposition 8.2]{Fuehring_bordism}, one shows that for $q > 0$, the composition
\[
 \Omega^{\SO}_{q-4} (\sB G) \stackrel{\psi({\rm pt.})_*}{\longrightarrow} \Omega^{\SO}_{q} \to \Omega^{\SO}_q/ \Tor_q
\]
is surjective.
Hence, the class $\xi \in \Omega_q^{\SO} / \Tor_q$ is represented by the total space of a fibre bundle with structure group $G$,
\[
 \C P^2 \hookrightarrow N^q \stackrel{\pi}{\longrightarrow} B^{q-4}.
 \]
Now, in \eqref{ses}, the class in $E^2_{p,q}$ represented by
\[
 [ f \colon M^p \to X] \cdot [N^q] = \big[M^p \times N^q \stackrel{\id \times \pi}{\longrightarrow} M^p \times B^{q-4} \stackrel{(x,y) \mapsto f(x)}{\longrightarrow} X\big] \in \mathscr{T}_n(X)_{(2)} \cap \mathscr{F}_{p,q}
\]
maps to $k ( h \otimes \xi) \in H_p\big(X; \Omega^{\SO}_q / \Tor_q \big)_{(2)}$, finishing the proof of Assertion~\ref{Assertion1}.

\begin{Assertion}\label{Assertion2}
 Each element in $H_p(X; \Tor_q) \subset E^2_{p,q}$ is represented by some class in $\mathscr{T}_n(X)_{(2)} \cap \mathscr{F}_{p,q}$.
\end{Assertion}

To show this, we use a homotopy-theoretic argument inspired by \cite{Stolz_splitting}.
For an abelian group~$A$, let ${\sf H}A$ denote the Eilenberg--MacLane spectrum for~$A$, and let $\mathscr{A}^*=H^*({\sf H}\F_2; \F_2)$ be the Steenrod algebra of stable $\F_2$-cohomology operations.
Recall that for any spectrum ${\sf M}$, the cohomology $H^*({\sf M}; \F_2)$ is an $\mathscr{A}^*$-module.

Let $\MSO$ be the oriented bordism spectrum and let $\MSO_{(2)}$ be its $2$-localisation.
Recall that the Pontrjagin--Thom construction yields a natural isomorphism
\[
 \Omega^{\SO}_*(X) \cong \pi_*(\MSO \wedge \Sigma^{\infty} X_+),
\]
where $\wedge$ denotes the smash product and $\Sigma^{\infty} X_+$ is the suspension spectrum of $X$ with a disjoint base point.
From this perspective, the homological orientation $u(X)_*$ from \eqref{homor} corresponds to a spectrum map
\begin{equation} \label{u}
 u \colon\ \MSO \to \mathsf{H}\Z
\end{equation}
and the map $\psi(X)_*$ from \eqref{psi} is identified with
\[
(T \wedge \id_{\Sigma^{\infty} X_+})_* \colon\ \pi_*\big( \Sigma^4 \MSO \wedge \Sigma^{\infty} \sB G_+ \wedge \Sigma^{\infty} X_+\big) \to \pi_*(\MSO \wedge \Sigma^{\infty}X_+) .
\]
Here, the spectrum map $T \colon \Sigma^4 \MSO \wedge \Sigma^{\infty} \sB G_+ \to \MSO$ is the composition
\[
 \Sigma^4 \MSO \wedge \Sigma^{\infty} \sB G_+ \stackrel{\id \wedge \alpha}{\longrightarrow} \Sigma^4 \MSO \wedge {\rm Th}(-T_v \sB H) \stackrel{\id \wedge \beta}{\longrightarrow} \Sigma^4 \MSO \wedge \Sigma^{-4} \MSO \stackrel{\gamma}{\longrightarrow} \MSO,
 \]
where ${\rm Th}(-T_v \sB H)$ is the Thom spectrum of the stable fiberwise normal bundle of \eqref{universal}, $\alpha \colon \Sigma^{\infty} \sB G_+ \to {\rm Th}(-T_v \sB H)$ is the Pontrjagin--Thom map, $\beta \colon {\rm Th}(-T_v \sB H) \to \Sigma^{-4} \MSO$ is the homological orientation, and $\gamma$ is multiplication in the ring spectrum $\MSO$.
For more information on the construction of $T$, compare \cite{Fuehring_bordism}, in particular, Section~3 and Lemma~6.1.

By \cite[Theorem 5]{Wall60}, also compare \cite[proof of~Theorem 14.1]{CF}, there are splittings of spectra
\begin{gather} \label{2localsplit}
 \MSO_{(2)} \simeq \sf{P} \vee \sf{Q} , \qquad \sf{Q} \simeq \mathsf{Q'} \vee \mathsf{H} \Z_{(2)},
 \end{gather}
where $\sf{P}$ is a sum of positive suspensions of ${\sf H}\F_2$ and $\mathsf{Q'}$ is a sum of $4\ell$-fold suspensions of ${\sf H}\Z_{(2)}$, $\ell > 0$.
The induced projection $\MSO_{(2)} \to \mathsf{H}\Z_{(2)}$ is the $2$-localisation of the map ${u \colon \MSO \to \mathsf{H}\Z}$ appearing in \eqref{u}.
In particular, denoting by \smash{$\widehat \MSO \subset \MSO$} the homotopy fibre of $u$, we have~\smash{$ ( \widehat \MSO)_{(2)} \simeq \sf{P} \vee \sf{Q'}$}.
Furthermore, with respect to \eqref{2localsplit}, we have
\begin{equation} \label{splitorbord}
 H^*(\MSO; \F_2) = H^*(\mathsf{P}; \F_2) \oplus H^*( \mathsf{Q}; \F_2) , \qquad \Omega^{\SO}_*(X)_{(2)} = \mathsf{P}_*(X) \oplus \mathsf{Q}_*(X) .
\end{equation}
For $X=\{\rm pt.\}$, the second equality specializes to \smash{$\big(\Omega^{\SO}_*\big)_{(2)} = \pi_*( \mathsf{P}) \oplus \pi_*(\mathsf{Q})$} with $\pi_*(\mathsf{P}) = \Tor_* \subset \Omega^{\SO}_*$.

Since the composition $u \circ T$ represents an element in $H^{-4}( \MSO \wedge \Sigma^{\infty} \sB G_+; \Z) = 0$, it is null-homotopic, and hence the map $T$ lifts to a map \smash{$\widehat T \colon \Sigma^4 \MSO \wedge \Sigma^{\infty} \sB G_+ \to \widehat \MSO$}.
According to \cite[Theorem 2.2]{Fuehring_bordism},\footnote{This reference works with the dual notion of $\mathscr{A}_*$-comodules where $\mathscr{A}_*$ is the dual $\F_2$-Steenrod algebra.} $\widehat T$ induces a split injection of $\mathscr{A}^*$-modules
\[
 \widehat T^* \colon\ H^*\big( \widehat\MSO; \F_2\big) \to H^{*}\big( \Sigma^4 \MSO \wedge \Sigma^{\infty} \sB G_+; \F_2\big) .
\]
Since \smash{$\big( \widehat \MSO\big)_{(2)} \simeq \sf{P} \vee \sf{Q'}$}, this implies that there exists a graded $\mathscr{A}^*$-module $M^*$ and a direct sum decomposition of $\mathscr{A}^{*}$-modules
\begin{equation} \label{splitting1}
 H^{*}\big(\Sigma^4 \MSO \wedge \Sigma^{\infty} \sB G_+; \F_2\big) \cong H^*(\mathsf{P}; \F_2) \oplus M^*,
\end{equation}
so that the restriction of $T^* \colon H^{*} (\MSO; \F_2) \to H^*\big(\Sigma^4 \MSO \wedge \Sigma^{\infty} \sB G_+; \F_2\big)$ to $H^*(\mathsf{P}; \F_2)$, see~\eqref{splitorbord}, is equal to the inclusion $H^*(\mathsf{P}; \F_2) \hookrightarrow H^{*}\big(\Sigma^4 \MSO \wedge \Sigma^{\infty} \sB G_+; \F_2\big)$.

Since $H^*(\mathsf{P}; \F_2)$ is a free $\mathscr{A}^*$-module, according to \cite[Theorem 2\,(a)]{Mar74}, which applies to the locally finite spectrum $\Sigma^4 \MSO \wedge \Sigma^{\infty} \sB G_+$, the splitting \eqref{splitting1} can be realised on the spectrum level, i.e., there is a spectrum $\mathsf{M}$ with $H^*(\mathsf{M} ;\F_2) \cong M^*$ as $\mathscr{A}^*$-modules and a homotopy equivalence (of unlocalized spectra)
\begin{equation} \label{Margolis1}
\Sigma^4 \MSO \wedge \Sigma^{\infty} \sB G_+ \simeq \mathsf{P} \vee \mathsf{M},
 \end{equation}
which induces the splitting \eqref{splitting1} in $\F_2$-cohomology.
Consider the map
\[
 \Xi \colon\ \mathsf{P} \stackrel{\eqref{Margolis1}}{\hookrightarrow} \Sigma^4 \MSO \wedge \Sigma^{\infty} \sB G_+ \stackrel{T}{\longrightarrow} \MSO \to \MSO_{(2)} \stackrel{\eqref{2localsplit}}{\simeq} \mathsf{P} \vee \mathsf{Q} .
 \]
In $\F_2$-cohomology, the composition of this map with the projection onto $\mathsf{P}$ induces the identity on $H^*(\mathsf{P}; \F_2)$, while the composition with the projection onto $\mathsf{Q}$ induces the $0$-map.
Hence, by \cite[Theorem 3]{Mar74} and using that $\mathsf{P}$ is a sum of suspensions of ${\sf H}\F_2$, the map $\Xi$ is homotopic to the inclusion $\mathsf{P} \to \mathsf{P} \vee \mathsf{Q}$ onto the first summand.

We conclude that the composition
\[
 \mathsf{P}_*(X) \to \pi_*\big( \Sigma^4 \MSO \wedge \Sigma^{\infty} \sB G_+ \wedge \Sigma^{\infty}X_+\big) \stackrel{\psi(X)_*}{\longrightarrow} \pi_*(\MSO \wedge \Sigma^{\infty} X_+) = \mathsf{P}_*(X) \oplus \mathsf{Q}_*(X)
\]
is equal to the inclusion $ \mathsf{P}_*(X) \subset \mathsf{P}_*(X) \oplus \mathsf{Q}_*(X)$.
In particular, in the decomposition \eqref{splitorbord}, we have $\mathsf{P}_*(X) \subset \mathscr{T}_*(X)_{(2)}$.

The splitting $\MSO_{(2)} \simeq \mathsf{P} \vee \mathsf{Q}$ in \eqref{2localsplit} induces a splitting $E^r_{p,q} = E^{\mathsf{P}, r}_{p,q} \oplus E^{\mathsf{Q}, r}_{p,q}$ of the Atiyah--Hirzebruch spectral sequence for $\Omega^{\SO}_*(X)_{(2)}$.
Clearly, the Atiyah--Hirzebruch spectral sequences for $\mathsf{P}$ and $\mathsf{Q}$ collapse at $E^2$, and we obtain an induced decomposition
\[
 H_p\big(X; \Omega^{\SO}_q\big)_{(2)} = E^2_{p,q} = E(\mathsf{P})^{2}_{p,q} \oplus E(\mathsf{Q})^{2}_{p,q} = H_p(X ; \pi_q( \mathsf{P})) \oplus H_p(X ; \pi_q( \mathsf{Q})),
\]
where
\[
 E(\mathsf{P})^{2}_{p,q} = \left( \mathscr{F}_{p,q} \cap \mathsf{P}_n(X)\right) / \left( \mathscr{F}_{p-1,q+1} \cap \mathsf{P}_n(X) \right).
\]
Since $H_p(X ; \pi_q( \mathsf{P})) = H_p(X; \Tor_q) \subset E^2_{p,q}$ and $\mathsf{P}_n(X) \subset \mathscr{T}_n(X)_{(2)}$, this finishes the proof of Assertion~\ref{Assertion2}.
\end{proof}

We obtain the following bordism theoretic description of $H_*(X ; \Z_{(2)})$ whose analogue for spin bordism and ${\rm ko}$-homology was shown in \cite[Theorem B\,(2)]{Stolz_splitting}.

\begin{Corollary} For each CW complex $X$, the orientation $u(X)_*$ induces an isomorphism
\[
 \left( \Omega^{\SO}_*(X) / \mathscr{T}_*(X) \right)_{(2)} \cong H_*(X; \Z)_{(2)} .
\]
\end{Corollary}

\begin{proof}[Proof of Theorem \ref{prop:nonspin}]
By the {\em bordism principle} \cite[Theorem 1.5]{R-S1995}, which follows from the {\em surgery principle} in scalar curvature geometry \cite{GL,SY}, and since the universal cover of $V$ is non-spin, it would suffice to show that the oriented bordism class
\[
 [ \varphi \colon V \to \sB \pi_1(V)] \in \Omega^{\SO}_*( \sB \pi_1(V))
\]
lies in $\Omega^{\SO, +}_*( \sB \pi_1(V))$.
We will achieve this goal for an odd multiple of $[ \varphi \colon V \to {\sf B} \pi_1(V)]$, while the rest of the proof relies on a version of the bordism principle for manifolds with Baas--Sullivan singularities worked out in \cite{Han20}.

By assumption, the class
\[
 \beta := [ \varphi \colon V \to \sB \pi_1(V) ] - [f \colon V' \to \sB \pi_1(V) ] \in \Omega_*( \sB \pi_1(V))
\]
is in the kernel of $u(\sB \pi_1 (V))_* \colon \Omega_*(\sB \pi_1(V)) \to H_*(\pi_1(V) ; \Z)$.
By Proposition \ref{crucial}, there is an odd $k \in \N$ with $k \beta \in \mathscr{T}_*(\sB \pi_1(V)) \subset \Omega^{\SO, +} ( \sB \pi_1(V)) $.
Write $k \beta = [ M \to \sB \pi_1(V)]$ where $M$ is~\yes.
Since $[f \colon V' \to \sB \pi_1(V) ] \in \Omega^{\SO,+}_*( \sB \pi_1(V))$ by assumption, we obtain
\[
 k\cdot [ \varphi \colon V \to \sB \pi_1(V) ] \in \Omega^{\SO, +}_*( \sB \pi_1(V)) .
\]

We now use the bordism principle for \yes-manifolds based on {\em manifolds with Baas--Sullivan singularities} \cite[Proposition 3.11]{Han20}.
Let $\Omega^{\SO, \mathscr{Q}}_*(-)$ denote oriented bordism with Baas--Sullivan singularities in $\mathscr{Q}$, where $\mathscr{Q} = (Q_0 = *, Q_1, Q_2, \ldots)$ is a family of closed oriented manifolds~$Q_i$ of dimension $4i$ which are equipped with fixed positive scalar curvature metrics and such that~${\Omega^{\SO}_* / \Tor = \Z [[Q_1], [Q_2], \ldots]}$.
Note that $\mathscr{Q}$ is a positive family of singularity types in the sense of \cite[Definition 3.2]{Han20}, also compare the remarks before \cite[Definition 3.12]{Han20}.

Recall from \cite{baas1} and \cite[Proposition 2.6]{Han20} that for any space $X$, there exists a canonical isomorphism
\[
 \Omega^{\SO, \mathscr{Q}}_*(X) \otimes \Z\big[\tfrac{1}{2}\big] \cong H_*\big(X; \Z\big[\tfrac{1}{2}\big] \big) .
 \]
This implies that the image of $\beta$ in $\Omega^{\SO, \mathscr{Q}}_*( \sB \pi_1(V)) \otimes \Z\big[\tfrac{1}{2}\big]$ is equal to zero.

Thus there exists some $d \geq 0$ such that the image of $2^d \cdot \beta$ in $\Omega^{\SO, \mathscr{Q}}_*( \sB \pi_1(V))$ is equal to $0$.
Let $s, t \in \Z$ be integers with $s \cdot k + t \cdot 2^d = 1$.
This implies
\[
 s k \cdot \beta = \beta \in \Omega^{\SO, \mathscr{Q}}_*( \sB \pi_1(V)) .
\]
Let $W \to \sB \pi_1(V)$ be a compact $\mathscr{Q}$-bordism between $s k \cdot \beta $ and $\beta$.
More precisely, $W$ is a compact oriented $\mathscr{Q}_n$-manifold for some $n \geq 0$ (see \cite[p.\ 503]{Han20}) and $W \to \sB \pi_1(V)$ compatible with the singularity structure of~$W$ (see \cite[Definition 2.3]{Han20}).
Furthermore, along the boundary $\partial_0 V$ of~$W$ with its induced orientation (see \cite[Definition 2.1]{Han20}), this map restricts to the disjoint union of $s$ copies of $M \to \sB \pi_1(V)$, of $\varphi \colon V^{\rm op}\to \sB \pi_1(V)$ and of $f \colon V' \to \sB \pi_1(V)$.
Here, $V^{\rm op}$ is the manifold $V$ with the opposite orientation.

By applying surgeries in the interior of $W$, we can assume that the inclusion $V \hookrightarrow W$ is a~$2$-equivalence, i.e., it induces bijections on $\pi_0$ and $\pi_1$ and a~surjection on $\pi_2$.
This follows by a~similar argument as in the smooth case, see \cite[Theorem 1.5]{R-S1995}.\footnote{For some details one may consult the proof of \cite[Proposition 3.1]{SchickZenobi}.}

Now, since $M$ and $V'$ are \yes, and $\mathscr{Q}$ is a positive family of singularity types, the $\mathscr{Q}$-bordism principle \cite[Proposition 3.11]{Han20} shows that the manifold $V$ is \yes, completing the proof of Theorem \ref{prop:nonspin}.
\end{proof}

\begin{Remark} Theorem~\ref{prop:nonspin} together with a short proof has already appeared in \cite[Theorem~4.11]{RS}.
However, this discussion is incomplete if $H_*(\pi_1(V); \Z)$ contains $2$-torsion, which may produce Tor terms in the K\"unneth formula for the $E^2$-terms \smash{$H_p\big( \sB \pi_1(V) ; \Omega^{\SO}_q\big)_{(2)}$} in the Atiyah--Hirzebruch spectral sequence for \smash{$\Omega^{\SO}_*( \sB \pi_1(V))_{(2)}$}.
This is covered by Assertion~\ref{Assertion2} in the proof of our Proposition \ref{crucial}.
\end{Remark}

\subsection*{Acknowledgements}

This paper is the result of conversations in New York in the fall of 2018 and 2022.
B.H.~acknowledges support from NYU, the IAS Princeton and the DFG-funded Special Priority Program 2026 {\em Geometry at Infinity}.
We are grateful to the referees for their many useful comments, which helped to improve the manuscript.

\pdfbookmark[1]{References}{ref}
\LastPageEnding

\end{document}